\documentclass{jag-l}

\usepackage{url}

\usepackage[utf8]{inputenc}

\usepackage[all]{xy}
\usepackage{pb-diagram, pb-xy}

\usepackage{amssymb}
\usepackage{amsthm}

\newtheorem{thm}{Theorem}[section]
\newtheorem{lem}[thm]{Lemma}
\newtheorem{cor}[thm]{Corollary}
\newtheorem{prop}[thm]{Proposition}
\newtheorem{rem}[thm]{Remark}
\newtheorem{ex}[thm]{Example}
\newtheorem{defn}[thm]{Definition}

\numberwithin{equation}{thm}

\newcommand{\bbC}{{\mathbb C}}
\newcommand{\bbG}{{\mathbb G}}

\newcommand{\bbN}{{\mathbb N}}

\newcommand{\bbQ}{{\mathbb Q}}
\newcommand{\bbR}{{\mathbb R}}

\newcommand{\bbZ}{{\mathbb Z}}

\newcommand{\sfH}{{\mathsf H}}

\newcommand{\frakg}{{\mathfrak g}}

\newcommand{\Qbar}{\overline{\mathbb Q}}

\newcommand{\Supp}{\mathrm{Supp}}

\newcommand{\calD}{{\mathcal D}}

\newcommand{\calH}{{\mathcal H}}

\newcommand{\calL}{{\mathcal L}}
\newcommand{\calM}{{\mathcal M}}
\newcommand{\calN}{{ N_\bfD}}
\newcommand{\calO}{{\mathcal O}}
\newcommand{\calS}{{\mathcal S}}

\newcommand{\Sp}{\mathrm{Sp}}

\newcommand{\Gl}{\mathrm{Gl}}

\newcommand{\Rep}{{\mathrm{Rep}}}
\newcommand{\Pic}{{\mathit{Pic}}}
\newcommand{\Alb}{{\mathit{Alb}}}
\newcommand{\Lie}{{\mathrm{Lie}}}

\newcommand{\Dbc}{{{D^{\hspace{0.01em}b}_{\hspace{-0.13em}c} \hspace{-0.05em} }}}

\newcommand{\Perv}{\mathrm{Perv}}

\newcommand{\ptau}{{^{p} \! \tau}}

\newcommand{\pH}{{^p \! H}}

\newcommand{\End}{{\mathit{End}}}
\newcommand{\Hom}{{\mathit{Hom}}}
\newcommand{\Ext}{{\mathit{Ext}}}
\newcommand{\id}{{\mathit{id}}}
\newcommand{\MHM}{{\mathrm{MHM}}}

\newcommand{\Stab}{{\mathit{Stab}}}

\newcommand{\one}{{\mathbf{1}}}

\newcommand{\bfA}{{\mathbf{A}}}

\newcommand{\bfD}{{\mathbf{D}}}
\newcommand{\bfG}{{\mathbf{G}}}

\newcommand{\bfN}{{\mathbf{N}}}
\newcommand{\bfP}{{\mathbf{P}}}

\newcommand{\bfS}{{\mathbf{S}}}
\newcommand{\bfT}{{\mathbf{T}}}
\newcommand{\bfV}{{\mathbf{V}}}

\newcommand{\overbar}[1]{\mkern 1.5mu\overline{\mkern-1.5mu#1\mkern-3mu}\mkern 3mu}

\newcommand{\bfDbar}{{\overbar{\mathbf{D}}}}

\newcommand{\bfPbar}{{\overbar{\mathbf{P}}}}

\newcommand{\kappabar}{{\mkern 3mu \overline{\mkern -3mu \kappa \mkern -1.5mu} \mkern 1.5mu}}
\newcommand{\etabar}{{\mkern 2.5mu \overline{\mkern -2.5mu \eta \mkern -2.5mu} \mkern 2.5mu}}
\newcommand{\sbar}{{\mkern 2mu \overline{\mkern -2mu  s } }}
\newcommand{\Sbar}{{\mkern 2mu \overline{\mkern -2mu S } }}
\newcommand{\Xbar}{{\mkern 2mu \overline{\mkern -2mu X \mkern -2mu} \mkern 2mu}}
\newcommand{\ibar}{{\overline i}}
\newcommand{\jbar}{{\overline{j}}}

\newcommand{\bfGl}{{\mathbf{Gl}}}

\newcommand{\CC}{{\mathit{CC}}}

\begin{document}

\title[Vanishing theorems for constructible sheaves]{Vanishing theorems for constructible sheaves on abelian varieties}
\author{Thomas Krämer}
\author{Rainer Weissauer}
\address{Mathematisches Institut\\ Ruprecht-Karls-Universit\"at Heidelberg\\ Im Neuenheimer Feld 288, D-69120 Heidelberg, Germany}
\email{tkraemer@mathi.uni-heidelberg.de}
\email{weissauer@mathi.uni-heidelberg.de}

\begin{abstract}
We show that the hypercohomology of most character twists of perverse sheaves on a complex abelian variety vanishes in all non-zero degrees. As a consequence we obtain a vanishing theorem for constructible sheaves and a relative vanishing theorem for a homomorphism between abelian varieties. Our proof relies on a Tannakian description for convolution products of perverse sheaves, and with future applications in mind we discuss the basic properties of the arising Tannaka groups.
\end{abstract}

\maketitle

\thispagestyle{empty}

\section{Introduction}

Let $X$ be a complex abelian variety, and denote by $\Dbc(X, \bbC)$ the derived category of bounded $\bbC$-sheaf complexes on $X$ with constructible cohomology sheaves (by Chow's theorem it makes no difference whether  for the notion of constructibility we use analytic or algebraic stratifications). By definition a complex $K\in \Dbc(X, \bbC)$ is semi-perverse if  its cohomology sheaves $\calH^{-i}(K)$ satisfy the estimate $\dim ( \Supp \, \calH^{-i}(K)) \leq i$ for all $i\in \bbZ$, and $K$ is called a perverse sheaf if both $K$ and its Verdier dual $DK$ are semi-perverse. Let
\[ \Perv(X,\bbC) \; \subset \; \Dbc(X, \bbC) \]
be the full subcategory of perverse sheaves. This is an abelian category, the core of the perverse $t$-structure on $\Dbc(X, \bbC)$ as defined in~\cite{BBD}. 

\medskip

The group structure on $X$ defines a convolution product on $\Dbc(X, \bbC)$ under which $\Dbc(X, \bbC)$ becomes a rigid symmetric monoidal triangulated category in a natural way, see \cite{WeBN} and~\cite{WeRem}. This convolution product does not preserve the full abelian subcategory of perverse sheaves, but we construct an abelian quotient category of $\Perv(X, \bbC)$ that is a Tannakian category in the sense of~\cite{DM} with respect to a tensor product induced by convolution. It turns out that the Tannakian property is essentially equivalent to a vanishing theorem for the hypercohomology $H^\bullet(X, P)$ of perverse sheaves $P\in \Perv(X, \bbC)$. To formulate this vanishing theorem, let $\Pi(X) = \Hom(\pi_1(X, 0), \bbC^*)$ denote the algebraic torus of characters of the fundamental group. Any $\chi \in \Pi(X)$ defines a local system $L_\chi$ of rank one on $X$, and we show

\medskip \begin{thm} \label{thm:P}
Let $P\in \Perv(X, \bbC)$. Then for all characters $\chi$ outside a finite union of translates of proper algebraic subtori of $\Pi(X)$ we have
\[ H^i(X, P\otimes_\bbC L_\chi)=0 
\quad \textnormal{\em for} \quad
i \neq 0.
\] 
\end{thm} \medskip

To make the statement of theorem~\ref{thm:P} more precise, let us introduce the following terminology. For abelian subvarieties $A \subseteq X$ let~$K(A) \subset \Pi(X)$ be the algebraic subtorus of all characters $\chi:\pi_1(X,0) \longrightarrow \bbC^*$ 
whose restriction to~$\pi_1(A,0)$ is trivial. By a {\em thin} set of characters we mean a finite union of translates $\chi_i \cdot K(A_i)$ for certain characters $\chi_i \in \Pi(X)$ and certain non-zero abelian subvarieties $A_i\subseteq X$. In these terms, we will show in section~\ref{sec:spectrum} that for any semisimple perverse sheaf $P$ the locus 
\[ \calS(P) \;=\; \bigl\{ \; \chi \in \Pi(X) \mid H^i(X, P\otimes_\bbC L_\chi) \neq 0 \;\; \textnormal{for some} \;\; i\neq 0  \;\bigr\} \]
is a thin subset of the character torus $\Pi(X)$. Writing $\calS(P)$ as a union of translates $\chi_i \cdot K(A_i)$ as above we will furthermore see that the $\chi_i$ can be chosen to be torsion characters,
if the perverse sheaf $P$ is of geometric origin in the sense of~\cite[6.2.4]{BBD}. In what follows, to save words we will say that a statement holds for {\em most} characters $\chi$ if it holds for all $\chi$ in the complement of a thin set of characters as defined above. 

\medskip

Theorem~\ref{thm:P} can easily be generalized to a relative vanishing theorem for a homomorphism of abelian varieties, see section~\ref{sec:relative-vanishing}.

\medskip

On algebraic tori, an analogue of theorem~\ref{thm:P} can be obtained from Artin's affine vanishing theorem and has been used in~\cite{GaL} for the construction of Tannakian categories of perverse sheaves. By way of contrast, for abelian varieties we define the Tannakian categories via a general construction of Andr\'e and Kahn~\cite{AK1}, which will allow to deduce theorem~\ref{thm:P} via the hard Lefschetz theorem and the theory of reductive (super)groups. Our proof in sections~\ref{sec:twists} -- \ref{sec:euler} is based on two ingredients. The first is a result of Deligne~\cite{DelCT} which characterizes rigid symmetric monoidal abelian categories and will be used to see that in the case at hand, the construction of Andr\'e and Kahn leads to a super Tannakian category in the sense of loc.~cit. To see that this category is in fact a Tannakian category in the usual sense, we require the second ingredient of the proof --- a classification of perverse sheaves with Euler characteristic zero in the spirit of~\cite{FK}, see proposition~\ref{prop:chi}. Here we use the theory of \mbox{$\calD$-modules}, and this is the only place where we need to work over the complex numbers. Except for section~\ref{sec:euler}, with the obvious modification of the notions {\em most} and {\em thin} our proof works in the same way for $\ell$-adic perverse sheaves on abelian varieties over the algebraic closure of a finite field as defined in~\cite{BBD}. Using this, our result has recently been generalized to the case of positive characteristic in~\cite{WeVanishing}.

\medskip

Via the Tannakian categories mentioned above, one can attach to any semisimple perverse sheaf $P\in \Perv(X, \bbC)$ a reductive complex algebraic group $G(P)$. In particular, for every smooth complex projective variety $Y$ with Albanese morphism $f: Y\longrightarrow X=\Alb(Y)$ we obtain a new invariant, the Tannaka group attached to the direct image complex $Rf_*(\bbC_Y[\dim(Y)])$ as in~\cite{WeConn}. Furthermore, the above groups are closely related to the moduli of abelian varieties~\cite{KrW}. Since therefore the Tannakian categories occuring in our proof are of independent interest, we explain in sections~\ref{sec:localization} through~\ref{sec:nearby} how their construction can be extended to the non-semisimple case, and we survey the basic properties of the arising Tannaka groups.

\medskip

Theorem~\ref{thm:P} can also be reformulated as a statement about constructible sheaves. Indeed, by d\'evissage with respect to the perverse $t$-structure and by Verdier duality one sees that theorem \ref{thm:P} is equivalent to the statement that any semi-perverse complex $K$ satisfies $H^i(X, K\otimes_\bbC L_\chi)=0$ for $i>0$ and most $\chi$. For any constructible sheaf $F$ the complex $K=F[\dim (\Supp \, F ) ]$ is semi-perverse, so we obtain

\medskip \begin{thm} \label{thm:F}
Let $F$ be a constructible sheaf of complex vector spaces on a complex abelian variety~$X$. Then for most characters $\chi$
we have
\[ H^i(X, F\otimes_\bbC L_\chi) \;=\; 0
\quad \textnormal{\em for} \quad
i \;>\; \dim (\Supp \, F).
\] 
\end{thm} \medskip

\noindent
This can be viewed as an analog of the Artin-Grothendieck affine vanishing theorem in the same way as one can consider the generic vanishing theorem of Green and Lazarsfeld \cite[th.~1]{GrL} as an analog of the Kodaira-Nakano vanishing theorem. Indeed the Green-Lazarsfeld theorem is a special case of our result as we will explain in more detail in section~\ref{sec:nakano}.

\section{A relative generic vanishing theorem} \label{sec:relative-vanishing}

Let $X$ be a complex abelian variety and $A\subseteq X$ an abelian subvariety with quotient $f: X\longrightarrow B = X/A$. Assuming theorem~\ref{thm:P} only on $A$, we obtain the following relative generic vanishing theorem; here the quantifier {\em most} can be read in the slightly stronger sense that it does not refer to the characters of $\pi_1(X, 0)$ but rather to their pull-back to the subgroup $\pi_1(A,0) \subseteq \pi_1(X, 0)$, see the remark preceding lemma~\ref{lem:RelativeSpectrum}. 

\medskip \begin{thm} \label{thm:relative}
Let $P$ be a perverse sheaf on $X$. Then for most $\chi$ the direct image complex $Rf_*(P\otimes_\bbC L_\chi)$ is a perverse sheaf on~$B$.
\end{thm} \medskip

{\em Proof.} Put $P_\chi = P\otimes_\bbC L_\chi$. By Verdier duality it will be enough to show that for most characters $\chi$ the direct image complex $Rf_*(P_\chi)$ satisfies the semi-perversity condition
\[ 
\dim \bigl( \Supp \; \calH^{-k}(Rf_*(P_\chi)) \bigr) \; \leq \; k \quad \textnormal{for all} \quad k\in \bbZ. 
\]
To check this condition, note that by lemma~2.4 and section 3.1 in~\cite{BF} we can find Whitney stratifications $X=\sqcup_\beta X_\beta$ and $B=\sqcup_\alpha B_\alpha$ such that 

\begin{itemize}
\item[a)]
the cohomology sheaves $\calH^{-i}(P_\chi) = \calH^{-i}(P) \otimes_\bbC L_\chi$ are locally constant on the strata $X_\beta$ for all $\beta$, $i$ and $\chi$, 
\item[b)] each~$f(X_\beta)$ is contained in some $B_\alpha$, and 
\item[c)] for all $\alpha, \beta$ with $f(X_\beta)\subseteq B_\alpha$ the restriction $f: X_\beta \to B_\alpha$ is smooth. 
\end{itemize}
By theorem 4.1 of loc.~cit.~then the restriction $\calH^{-k}(Rf_*(P_\chi))|_{B_\alpha}$ is locally constant for all $\alpha$, $k$ and $\chi$. Since there are only finitely many strata $B_\alpha$ and since $\calH^{-k}(Rf_*(P_\chi))\neq 0$ for only finitely many $k$, it follows that if the direct image complex $Rf_* (P_\chi)$ were not semi-perverse for most $\chi$, then we could find $\alpha$ and~$k$ such that
\begin{itemize}
\item[d)] $\dim(B_\alpha) > k$ (where as usual by the dimension of a constructible subset we mean the maximum of the dimensions of the irreducible components of its closure), and
\item[e)] $\calH^{-k}(Rf_*(P_\chi))_b \;\neq \; 0$ for all points $b\in B_\alpha(\bbC)$ and all $\chi$ in a set of characters which is not thin in the sense of the introduction. 
\end{itemize}
Indeed, if a property does not hold for most characters, then by definition it fails on a set of characters which is not thin. Fixing $\alpha$ and $k$ as above, we now argue by contradiction.

\medskip

Fix $b\in B_\alpha(\bbC)$. Consider the fibre $F_b = f^{-1}(b)$, and for arbitrary $\chi$ denote by $M_b = P_\chi|_{F_b}$ the restriction of $P_\chi$ to $F_b$ (we suppress the character twist in this notation). For the perverse cohomology sheaves
\[ M_b^r \;=\; \pH^{-r}(M_b)  \]
we have the spectral sequence
\[
 E_2^{rs} \;=\; H^{-s}(F_b, M_b^r) \; \Longrightarrow \; H^{-(r+s)}(F_b, M_b) \;=\; \calH^{-(r+s)} (Rf_* (P_\chi))_b.
\]
Theorem~\ref{thm:P} for $F_b\cong A$ shows that for most $\chi$ we have $H^{-s}(F_b, M_b^r) = 0$ for all $s\neq 0$ and all $r\in \bbZ$. For such $\chi$ the spectral sequence degenerates, i.e.
\[
 \calH^{-k} (Rf_* (P_\chi))_b \;=\; H^0(F_b, M_b^k).
\]
On the other hand by e) we can assume $\calH^{-k} (Rf_* P_\chi)_b \neq 0$. By the above then $M_b^k \neq 0$. Since $M_b^k = \pH^0(M_b[-k])$, it follows by definition of the perverse $t$-structure that
\[
 \dim \bigl( \Supp \; \calH^{-i}(M_b) \bigr) \; = \; i-k \;\geq \; 0 \quad \textnormal{for some} \quad i\in \bbZ. 
\]
Now by a) the support of $\calH^{-i}(P_\chi)$ is a union of certain strata $X_\beta$, so using the above dimension estimate and the definition of $M_b = P_\chi|_{F_b}$ we find a stratum $X_\beta\subseteq \Supp\; \calH^{-i}(P_\chi)$ with $\dim(F_b\cap X_\beta) = i-k$. Since by b) and c) the stratum $X_\beta$ is equidimensional over $B_\alpha$, it follows that 
\[
 \dim \bigl( \Supp \; \calH^{-i}(P_\chi) \bigr) \;\geq \; \dim(X_\beta) \; = \; i-k+\dim(B_\alpha).
\]
But $\dim(B_\alpha) > k$ by property d), so it follows that the perverse sheaf $P_\chi$ is not semi-perverse, a contradiction.
\qed

\medskip

Note that in the proof of theorem~\ref{thm:relative} we have only used theorem~\ref{thm:P} for the fibres $f^{-1}(b) \cong A$ but not for $X$ itself. In fact, using this observation and assuming theorem~\ref{thm:P} only for simple abelian varieties, one can by induction on the dimension deduce for arbitrary abelian varieties a weaker version of theorem~\ref{thm:P} where the quantifier {\em most} is replaced by {\em generic} \cite{WeGauss}.

\section{Kodaira-Nakano-type vanishing theorems} \label{sec:nakano}

From theorem~\ref{thm:P} one easily recovers stronger versions of the vanishing theorems of Green and Lazarsfeld as follows. Let~$Y$ be a compact connected K\"ahler manifold of dimension $d$ whose Albanese variety $\Alb(Y)$ is algebraic, and denote by 
\[ f: \; \; Y \; \longrightarrow \; X \; = \; \mathrm{Alb}(Y) \] 
the Albanese morphism. To pass from coherent sheaves to constructible sheaves, recall that every coherent line bundle $\calL \in \Pic^0(X)$ admits a flat connection. The horizontal sections for any such connection form a local system $L_\chi$ where $\chi: \pi_1(X,0) \longrightarrow \bbC^*$ is a character with $\calL \cong  L_\chi  \otimes_\bbC  \calO_X$.  

\medskip

For a given line bundle $\calL \in \Pic^0(X)$, the set of all characters $\chi$ with the above property is a torsor under the group $H^0(X, \Omega^1_X)$. Indeed, this follows from the truncated exact cohomology sequence
\[
 0 \; \longrightarrow \; H^0(X, \Omega_X^1) \; \longrightarrow \; H^1(X, \bbC^*) \; \longrightarrow \; \Pic^0(X) \; \longrightarrow \; 0
\]
attached to the exact sequence $0\rightarrow \bbC_X^* \rightarrow \calO_X^* \rightarrow \Omega_{X, \mathit{cl}}^1 \rightarrow 0$ where $\Omega_{X, \mathit{cl}}^1$ denotes the sheaf of closed holomorphic $1$-forms. On the other hand, from the point of view of Hodge theory it is better to restrict our attention to unitary characters $\chi: \pi_1(X, 0) \longrightarrow U_1 = \{ z\in \bbC^* \mid |z|=1\}$, which has the extra benefit that it makes the passage from coherent to constructible sheaves unique: Comparing the exponential sequences $0\rightarrow \bbZ_X \rightarrow \bbR_X \rightarrow U_{1,X} \rightarrow 0$ and $0\rightarrow \bbZ_X \rightarrow \calO_X \rightarrow \calO_X^* \rightarrow 0$ one sees that the morphism
\[
 H^1(X, U_1) \; \stackrel{\cong}{\longrightarrow} \; \Pic^0(X)
\]
is an isomorphism, so for every line bundle $\calL \in \Pic^0(X)$ there is a unique unitary character $\chi$ with $\calL \cong L_\chi \otimes_\bbC \calO_X$. Concerning the applicability of theorem~\ref{thm:P} in this unitary context, we remark that the intersection of any thin subset of $\Pi(X)$ with the set of unitary characters is mapped via the above isomorphism to a thin subset of $\Pic^0(X)$, with the definition of {\em thin} and~{\em most} being extended in the obvious way to the Picard group.


\medskip

In what follows we put $X_n = \{ x\in X \mid \dim (f^{-1}(x)) = n \}$ for $n \in \bbN_0$ and consider the integer
\[
  w(Y) \;=\; \min \bigl\{ 2d - (\dim(X_n)+2n) \mid n \in \bbN_0, \, X_n \neq \varnothing \bigr\}.
\]
Notice that~$w(Y) \leq d$. Indeed, for some $n$ the preimage $f^{-1}(X_n)$ is dense in~$Y$ so that $d=\dim (f^{-1}(X_n))=\dim(X_n)+n$, hence $2d - (\dim(X_n)+2n)$ is equal to $2d - (d+n) = d- n \leq d$ as required. 

\medskip
In particular, the morphism $f$ is semi-small in the sense of~\cite[III.7]{KW} if and only if $w(Y) = d$. Furthermore, for local systems $E$ on $Y$ one sees as in loc.~cit.~that the complex
\[ Rf_* (E[2d-w(Y)]) \quad \textnormal{is semi-perverse.} \]
Hence theorem~\ref{thm:P} implies the following version of the Kodaira-Nakano type vanishing theorem of Green and Lazarsfeld~\cite[th.~2]{GrL}.

\medskip 

\begin{thm} \label{thm:nakano}
Let $E$ be a unitary local system on $Y$. Then for most $\calL$ in~$Pic^0(Y)$ we have
\[ 
 H^p(Y, \Omega_Y^q(E\otimes_\bbC \calL)) = 0
 \quad
 \textnormal{\em for} 
 \quad
 p+q < w(Y).
\]
\end{thm} \medskip

{\em Proof.} The morphism $f^*: \Pic^0(X) \longrightarrow \Pic^0 (Y)$ is an isomorphism by construction of the Albanese variety~\cite[p.~553]{GH}, so every $\calL \in \Pic^0(Y)$ arises as the pull-back of some $\calM \in \Pic^0(X)$. As explained above, there is a unique unitary character $\chi$ such that
\[ \calM \; \cong \; \calO_X\otimes_\bbC L_\chi. \] 
Then $\calL \cong f^*(\calM ) \cong \calO_Y\otimes_\bbC f^*(L_\chi)$. Since all the occuring local systems are unitary, Hodge theory implies that
\medskip
\[
 \bigoplus_{p+q=k} H^p(Y, \Omega_Y^q(E\otimes_\bbC \calL)) \; \cong \; H^k(Y, E\otimes_\bbC f^* (L_\chi)).
\]
Putting $K=Rf_*E[2d-w(Y)]$ we can identify the cohomology group on the right hand side with the group
$
  H^{k-2d+w(Y)}(X, K_\chi).
$
Since the direct image complex $K_\chi$ is semi-perverse, theorem \ref{thm:P} shows that for $k > 2d-w(Y)$ and most characters $\chi$ the above group vanishes. The theorem now follows by an application of Serre duality.\qed

\medskip

For a similar result in this direction, let us consider for $n\in \bbN_0$ the closed analytic subsets 
\[ \overline{X}_n = \{x\in X \mid \dim (f^{-1}(x)) \geq n\} \quad \textnormal{and} \quad \overline{Y}_n = f^{-1}(\overline{X}_n), \] 
and put $d_n = \dim(\overline{Y}_n)$ with the convention that $d_n=-\infty$ for $\overline{Y}_n =\emptyset$. Then our vanishing theorem implies the following

\medskip \begin{thm} \label{thm:level} Suppose that $p+q= d-n$ for some $n\geq 1$. Then for most line bundles $\calL$ in $\Pic^0(Y)$ we have
\[ H^p(Y, \Omega_Y^q(\calL)) \ =\ 0 \quad \textnormal{\em unless} \quad d-d_n \leq\ p,q \ \leq d_n-n. \]
\end{thm} \medskip

{\em Proof.} By Serre duality the claim of the theorem is equivalent to the statement that if~$p+q=d+n$ for some $n\geq 1$, then $H^p(Y, \Omega_Y^q(\calL))=0$ for most $\calL$ unless the Hodge types satisfy the estimates 
\[ d+n-d_n \; \leq \; p, q \; \leq \; d_n. \]
In fact it will suffice to establish the upper estimate $p,q\leq d_n$. Since we have $p+q=d+n$ by assumption, the lower estimate is then automatic.

\medskip
 
The decomposition theorem for compact K\"ahler manifolds~\cite[th.~0.6]{SaDT} says that $Rf_*\bbC_Y[d] \cong \bigoplus_m M_m[-m]$ where each $M_m$ is a pure Hodge module on $X$ of weight $m+d$ in the sense of~\cite{Sa}. Furthermore, for any unitary character $\chi$ with complex conjugate $\bar{\chi}$ the local system $L_\chi \oplus L_{\bar{\chi}}$ of rank two has an underlying real structure and hence can be viewed as a real Hodge module of weight zero in a natural way. So for any real Hodge module $M$ on $X$ also $M_{\chi, \bar{\chi}} = M_\chi \oplus M_{\bar{\chi}}$ is a real Hodge module. This being said, by theorem~\ref{thm:P} we have
$$H^{d+n}(Y,f^*(L_{\chi} \oplus L_{\bar{\chi}})) \;\cong\; H^{n}(X,(Rf_*\bbC_Y[d])_{\chi, \bar{\chi}}) \;\cong\; H^0(X,(M_{n})_{\chi, \bar{\chi}})$$
for most unitary characters $\chi$. The formalism of Hodge modules equips the cohomology group on the right hand side with a pure $\bbR$-Hodge structure of weight $n+d$ compatible with the natural one on the left hand side. We are looking for bounds on the types $(p, q)$ in this Hodge structure. 

\medskip

One easily checks that  $\Supp(M_n)\subseteq \overline X_n$, so $M_{n}[-n]$ is a direct summand of $Rf_* \bbC_{\bar{Y}_n}[d]$ by base change. To control the Hodge structure on twists of the cohomology of this direct image, let $\pi:\tilde Y\to Y$ be a composition of blow-ups in smooth centers that gives rise to an embedded resolution of singularities $\tilde Y_n = \pi^{-1}(\overline Y_n) \to \overline Y_n$, see \cite{Hi} or \cite[th.~10.7]{BM}. Then $\bbC_Y[d]$ occurs as a direct summand of the complex $R\pi_*\bbC_{\tilde Y}[d]$ by the decomposition theorem, so the restriction $\bbC_{\overline{Y}_n}[d]$ is a direct summand of $R\pi_* \bbC_{\tilde Y_n}[d]$. It then follows that $M_{n}[-n]$ is a direct summand of $Rf_* R\pi_* \bbC_{\tilde Y_n}[d]$, and we get an embedding
\[ H^0(X,(M_{n})_{\chi, \bar{\chi}}) \; \hookrightarrow \; H^{d+n}(\tilde Y_n,\pi^*f^*(L_{\chi} \oplus L_{\chi^{-1}})).
\]
But the Hodge types $(p, q)$ on the right hand side satisfy $p,q \leq \dim(\tilde{Y}_n) = d_n$ as one may check from the Hodge theory of compact K\"ahler manifolds with coefficients in unitary local systems.
\qed

\medskip

The above result contains the generic vanishing theorem of Green and Lazarsfeld~\cite[second part of th.~1]{GrL} as the special case $q=0$. Indeed, for any $p<\dim(f(Y))$ the number $n=d-p$ is larger than the dimension of the generic fibre of the Albanese morphism, hence $d_n<d$ so that $H^p(Y, \calL)=0$ for most $\calL$ by theorem~\ref{thm:level}. If $Y$ is algebraic, the theorem also holds more generally for $H^p(Y, \Omega_Y^q(E\otimes_\bbC \calL))$ with a unitary local system $E$ on $Y$.


\medskip

In general the bounds in the above theorem are strict: If~$d=4$ and if~$Y$ is the blow-up of $X$ along a smooth algebraic curve $C\subset X$ of genus~$\geq 2$, then one has $w(Y)=d_1=3$ but $H^2(Y, \Omega_Y^1(\calL))\neq 0$ for all non-trivial line bundles~$\calL$ as explained in \cite[top of p.~402]{GrL}.

\section{Character twists and convolution} \label{sec:twists}

We now introduce the notions of character twists and convolution, and we show that the two are compatible with each other. This will play a crucial role for our proof of theorem~\ref{thm:P} and for the construction of the Tannakian categories mentioned in the introduction. Indeed, the tensor product in these Tannakian categories will be given by the convolution product, but the fibre functors on them will only be constructed after a general character twist. 

\medskip

For the rest of this paper we work in the following setting. Let $X$ be an abelian variety over an algebraically closed field $k$ which has characteristic zero or is the algebraic closure of a finite field.  As in \cite{BBD} we consider the derived category~$\Dbc(X, \Lambda)$ of bounded complexes of $\Lambda$-sheaves on $X$ with constructible cohomology sheaves, where $\Lambda$ is either a subfield of~$\Qbar_l$ for some fixed prime number $l\neq char(k)$ or a subfield of $\bbC$, if we are working over the base field $k=\bbC$. We will denote by $\pi_1(X, 0)$ the \'etale fundamental group in the former and the topological fundamental group in the latter case. In the \'etale setting, by a character $\chi: \pi_1(X, 0) \longrightarrow \Lambda^*$ we always mean a continuous character whose image is contained in a finite extension field  of~$\bbQ_l$. Any such character defines a local system $L_\chi$ of rank one, and for $K\in \Dbc(X, \Lambda)$ we consider the corresponding character twist $K_\chi = K\otimes_\Lambda L_\chi$.
 
\medskip

Let $a: X\times X \longrightarrow X$ be the group law. Then $\Dbc(X, \Lambda)$ is a $\Lambda$-linear rigid symmetric monoidal category with respect to the convolution product
\[
 *: \;\; 
 \Dbc(X, \Lambda) \times \Dbc(X, \Lambda) \; \longrightarrow \; \Dbc(X, \Lambda), \quad
 K_1*K_2 = Ra_*(K_1\boxtimes K_2),
\]
see \cite[sect.~2.1]{WeBN} and \cite{WeRem}. The adjoint dual of an object $K$ in $\Dbc(X, \Lambda)$ is given in terms of its Verdier dual $DK$ by
\[
 K^\vee \;=\; (-\id_X)^* DK,
\]
and the unit object $\one$ of $\Dbc(X, \Lambda)$ is the skyscraper sheaf $\delta_0$ of rank one with support in the origin. Every skyscraper sheaf $K=\delta_x$ of rank one, supported in a point $x\in X(\bbC)$, is an invertible object in the sense that the evaluation morphism $K^\vee * K \longrightarrow \one$ is an isomorphism. Over the base field~$k=\bbC$ every invertible object has this form, as we will see in proposition~\ref{prop:chi}(b). 

\medskip

 To stress the symmetric monoidal structure on $\Dbc(X, \Lambda)$, we will sometimes use the notation $(\Dbc(X, \Lambda), *)$. We claim that twisting by a character defines on this symmetric monoidal category a tensor functor ACU in the sense of \cite[sect.~I.4.2.4]{Ri}. This fact will be crucial later on, though its proof is formal and may be skipped at a first reading.

\medskip \begin{prop} \label{prop:tensorfunctor}
For any character $\chi$, the autoequivalence $K\mapsto K_\chi$ of the category $\Dbc(X, \Lambda)$ defines a tensor functor ACU compatible with degree shifts and perverse truncations.
\end{prop} \medskip

{\em Proof.} The functor $K\mapsto K_\chi=K\otimes_\Lambda L_\chi$ preserves semi-perversity, so it is $t$-exact with respect to the perverse $t$-structure since $D(K_\chi)\cong D(K)_{\chi^{-1}}$. It remains to check tensor functoriality. Clearly $\one_\chi = \one$.

\medskip

Depending on the context,  put $R=\bbZ_l$, $R=\bbZ$ or $R=\bbZ/n\bbZ$ (the case where $p=char(k)$ divides  $n$ is included). The group law $a: X\times X\to X$ induces on cohomology the diagonal map 
\[
a^*: \; H^1(X, R) \to H^1(X\times X, R) = H^1(X, R)\oplus H^1(X, R), \; x \mapsto (x, x).
\]
In the first two cases use the formula preceding lemma~15.2 in~\cite{Mi}. In the last case notice that $\bbZ/n\bbZ \cong \mu_n$ for $(n,p)=1$  since $k$ is algebraically closed, and  $H^1(X,\mu_n) \cong Pic^0(X)[n]$ by \cite[cor.~III.4.18]{MiEC}. Thus for $(n,p)=1$ the claim follows since $a^*({\calL})\cong pr_1^*({\calL}) \otimes pr_2^*({\calL})$ holds for line bundles ${\calL}\in Pic^0(X)$, see \cite[prop.~9.2]{Mi}. 
On the other hand, 
$H^1(X,\bbZ/n\bbZ) \cong H^1(X,W_r)^{F}$ for $n=p^r$ by~\cite[prop.~13]{S1}. In this case,
the result follows by taking Frobenius invariants in $H^1(X\times X,W_r)= H^1(X,W_r) \oplus H^1(X,W_r)$, see \cite[p. 136]{S2}. 

\medskip

Now we have
$H^1(X, R) = \Hom(\pi_1(X, 0), R)$,
where in the \'etale setting we require the homomorphisms to be continuous; see \cite[rem.~15.5]{Mi} for $R=\bbZ_l$ and \cite[p. 50]{S1} for $R=\bbZ/n\bbZ$. If we write the group structure on fundamental groups additively, it follows that
\[ a_*: \pi_1(X, 0) \times \pi_1(X, 0) = \pi_1(X\times X, 0)\to \pi_1(X,0) \]
is the addition morphism $(x, y) \mapsto x+y$. For $\psi \in \Hom(\pi_1(X, 0), R)$ this implies 
$\psi(a_*(x, y)) = \psi(x+y) = \psi(x) + \psi(y)$, i.e.~$\psi \circ a_* = \psi \boxtimes \psi$
as an additive character on $\pi_1(X, 0)\times \pi_1(X, 0) = \pi_1(X\times X, 0)$. 
For multiplicative characters $\chi: \pi_1(X,0) \to \Lambda^*$ this implies
\[
 \chi(a_*(x, y)) = \chi(x+y) = \chi(x) \cdot \chi(y), \quad \textnormal{i.e.} \quad
 \chi \circ a_* \;=\; \chi \boxtimes \chi.
\]
Indeed, for $\Lambda \subseteq \bbC$  one has $\Hom(\pi_1(X, 0), R) \otimes_R \bbC^* = \Hom(\pi_1(X, 0), \bbC^*)$ taking $R=\bbZ$. For $\Lambda \subseteq \Qbar_l$ any multiplicative character $\chi$ takes values in the unit group $E^*\cong {\bbZ} \times { F}^* \times U$, where $F$ is the residue field of a finite extension field $E$ of $\bbQ_l$ and $U$ is its group of $1$-units. By continuity, $\chi=\chi_{F}\cdot \chi_U$ for characters $\chi_{F}$ and $\chi_U$ with values in ${F}^*$ resp.~$U$. The character $\chi_U$ can be handled as above, and the discussion for the character $\chi_{F}$ is covered by the case $R=\bbZ/n\bbZ$ with $n=|{F}^*|$.

\medskip

For the local system $L=L_\chi$ defined by a character $\chi: \pi_1(X, 0) \to \Lambda^*$ this gives an isomorphism on $X\times X$
\[ \varphi: \; a^*L  \stackrel{\sim}{\longrightarrow}  L \boxtimes L. \] 
Note that $\varphi$ is uniquely determined up to multiplication by an element of~$\Lambda^*$. In what follows,  we fix a choice of $\varphi$ once and for all. The choice of $\varphi$ will not matter for the commutativity of the diagrams below, as long as we use the same $\varphi$ consistently. 
However, since a tensor functor is not determined by the underlying functor alone, different choices of $\varphi$ give different (but isomorphic) tensor functors.
For us, it is most convenient to fix a trivialization $\lambda: L_0\cong \Lambda$ of the stalk $L_0$ at the origin $0$ of $X$, and to require that the stalk morphism
$  \varphi_{0}: a^*L_{(0,0)} \to    (L \boxtimes L)_{(0,0)} = L_0 \otimes_\Lambda L_0 $
at the origin $(0,0)$ of $X\times X$ makes the following diagramm commutative:
$$    \xymatrix{ (a^* L)_{(0,0)}\ar@{=}[d]  \ar[r]^{\varphi_0} &    L_0 \otimes_\Lambda L_0 \ar[d]^{\lambda \otimes \lambda}\cr         
 L_0  \ar[r]^{\lambda} &  \Lambda   \cr} $$
Here $L_0 = e_X^*(L) = e_{X^2}^*a^*(L)=(a^*L)_{(0,0)}$ since $a\circ e_{X^2} = e_X$ holds for the 
unit sections $e_X: \{0\}\to X$ and $e_{X^2}: \{(0,0)\} \to X^2$. For the unique $v\in L_0$ such that $\lambda(v)=1$, we have $\varphi_0^{-1}(\alpha\cdot v \otimes \beta\cdot v)= \alpha\beta \cdot v$ for $\alpha,\beta\in \Lambda$.

\medskip

Let $A, B\in \Dbc(X, \Lambda)$, and let $p_1, p_2: X\times X \to X$ be the projections onto the two factors. Using our fixed choice of $\varphi$, we get an  isomorphism
\[
 \psi: \; (A*B)_\chi  \stackrel{\sim}{\longrightarrow}  A_\chi*B_\chi
\]
defined by the commutative diagram
\[
\xymatrix@C=5em{
 (A*B)_\chi \ar@{=}[d] \ar@{..>}[r]^-{\psi} & A_\chi * B_\chi \\
 (Ra_*(A\boxtimes B)) \otimes L \ar@{=}[d] & Ra_*((A\otimes L)\boxtimes (B\otimes L)) \ar@{=}[u] \\ 
 Ra_*((A\boxtimes B)\otimes a^* L) \ar[r]^-{Ra_*(\id\, \otimes \varphi)} & Ra_*((A\boxtimes B)\otimes (L\boxtimes L)) \ar[u]^-\cong_-{Ra_*(\id \, \otimes S' \otimes \id)}
}
\]
where by $S': p_2^* (B) \otimes p_1^* (L) \stackrel{\sim}{\longrightarrow} p_1^* (L) \otimes p_2^* (B)$ we denote the symmetry constraint of the tensor product.

\medskip

The isomorphisms $\psi$ are compatible with the symmetry constraint~$S$ of the symmetric monoidal category $(\Dbc(X, \Lambda), *)$, i.e.~for all $A, B$ in~$\Dbc(X, \Lambda)$ the diagram
\[
\small
\xymatrix@C=5em{
 (A*B)_\chi \ar[d]_{S_\chi} \ar[r]^{\psi} & A_\chi*B_\chi \ar[d]^{S} \\
 (B*A)_\chi \ar[r]^{\psi} & B_\chi*A_\chi
}
\]
is commutative. Indeed, unravelling the definitions, the commutativity of the above diagram is equivalent to the commutativity of the diagram
\[
\xymatrix@C=3em{
 a^* L \ar[r]^-\varphi \ar@{=}[d] & L\boxtimes L \ar@{=}[r] & p_1^* L \otimes p_2^* L  \\
 \sigma^* a^* L \ar[r]^-{\sigma^*(\varphi)} & \sigma^* (L\boxtimes L) \ar@{=}[r] & p_2^* L \otimes p_1^* L \ar[u]_-{S'}^-\cong
}
\]
where $\sigma: X\times X \to X\times X$ is the morphism $(x,y)\mapsto (y, x)$ and $S'$ is the symmetry constraint of the tensor product. Since our diagram commutes up to a scalar in $\Lambda^*$, it suffices to check commutativity on the stalks at $(0,0)$. This boils down to the property $(\lambda\otimes \lambda)(u\otimes v) = (\lambda\otimes \lambda)(v\otimes u)$.

\medskip

The isomorphisms $\psi$ are also compatible with the associativity constraint of the symmetric monoidal category $(\Dbc(X, \Lambda), *)$. Indeed, we know by strictness \cite[p.~11]{WeBN} that the associativity constraints are the identity morphisms, so it suffices that  the diagram
\[
\small
\xymatrix@C=3em{
 ((A*B)*C)_\chi \ar@{=}[d] \ar[r]^-\psi
 & ((A*B)_\chi)*C_\chi \ar[r]^-{\psi * \id}
 & (A_\chi*B_\chi)*C_\chi \ar@{=}[d] \\
 (A*(B*C))_\chi \ar[r]^-\psi
 & A_\chi * ((B*C)_\chi) \ar[r]^-{\id \, * \psi}
 & A_\chi*(B_\chi*C_\chi) \\
}
\]
commutes for all $A, B, C\in \Dbc(X, \Lambda)$. Writing
\[ ((A*B)*C)_\chi = Ra_* R(a\times \id)_* (((A\boxtimes B)\boxtimes C) \otimes (a\times \id)^* a^* L) \]
and similarly for the other convolutions, the commutativity of the diagram becomes equivalent to the commutativity of the diagram on $X\times X\times X$ 
\[
\small
\xymatrix@C=3em@M=0.5em{
 (a\times \id)^* a^* L \ar@{=}[d] \ar[r]^-{(a\times \id)^* \varphi}
 & (a\times \id)^* L\boxtimes L =  a^* L \boxtimes L \ar[r]^-{\varphi \boxtimes \id}
 & (L\boxtimes L) \boxtimes L \ar@{=}[d] \\
 (\id \times a)^* a^* L \ar[r]^-{(\id \times a)^* \varphi}
 & (\id \times a)^* L\boxtimes L = L\boxtimes a^* L \ar[r]^-{\id\, \boxtimes \varphi}
 & L\boxtimes (L\boxtimes L) 
}
\]
Again it suffices to check the commutativity on stalks at $(0,0,0)$.
The upper arrow becomes the composition 
\[ (\varphi \otimes id)\circ \varphi: \quad L_0 \to L_0\otimes_\Lambda L_0
\to (L_0\otimes_\Lambda L_0)\otimes_\Lambda L_0. \]
Its inverse maps $(\alpha\cdot v \otimes \beta \cdot v) \otimes
\gamma\cdot v $ to $(\alpha\beta)\gamma \cdot v$. By a similar
 computation for the lower row, the commutativity of the diagram hence boils down to the associativity law $(\alpha\beta)\gamma = \alpha(\beta\gamma)$ of the field $\Lambda$.
\qed

\medskip

As a by-product, the tensor functoriality provides a simple proof of the following result from~\cite{IlEuler}.

\medskip \begin{cor} \label{cor:chi}
For $K\in \Dbc(X, \Lambda)$ the Euler characteristic of $K_\chi$ does not depend on the character $\chi$.
\end{cor} \medskip

{\em Proof.}
In~\cite[lemma~8 on p.~28]{WeBN} it has been deduced from the K\"unneth formula that hypercohomology defines a tensor functor ACU
\[
 H^\bullet(X, -): \;\; (\Dbc(X, \Lambda), *) \; \longrightarrow \; (Vect_\Lambda^s, \otimes^s)
\]
where the right hand side denotes the rigid symmetric monoidal category of super vector spaces over $\Lambda$, i.e.~the category of $\bbZ/2\bbZ$-graded vector spaces where the symmetry constraint is twisted by the usual sign rule. Hence the Euler characteristic of $K$ is equal, as an element of $\End_{\Dbc(X, \Lambda)}(\one) = \Lambda$, to the composite morphism
\vspace*{-0.5em}
\[
\xymatrix@C=3em{
 \one \ar[r]^-{coev_K} &
  K*K^\vee \ar[r]^-{S_{K, K^\vee}} &
  K^\vee * K \ar[r]^-{ev_K} &
 \one,
}
\]
and as such it is invariant under character twists by proposition~\ref{prop:tensorfunctor}. Here we denote by $ev_K$, $coev_K$ and $S_{K,K^\vee}$ the evaluation, coevaluation and symmetry constraint in the rigid symmetric monoidal category $(\Dbc(X, \Lambda), *)$. \qed

\medskip

\section{An axiomatic framework} \label{sec:setting}

Since the Tannakian constructions to be given below are of interest also in more general situations than in the proof of theorem~\ref{thm:P}, for the rest of this chapter we work in the following axiomatic setting. Let $(\bfD, *)$ be a $\Lambda$-linear rigid symmetric monoidal category with unit object $\one$, and let
\[
 rat: \; \; (\bfD, *) \; \longrightarrow \; (\Dbc(X, \Lambda), *)
\]
be a faithful $\Lambda$-linear tensor functor ACU. The notation $rat$ is motivated by the case where $k=\bbC$, $\Lambda = \bbQ$ and where $\bfD=D^b(\MHM(X))$ is the bounded derived category of the category $\MHM(X)$ of mixed Hodge modules~\cite{Sa}.

\medskip

For $K\in \bfD$ we denote by $H^\bullet(X, K)$ and by $\chi(K)$ the hypercohomology resp.~the Euler characteristic of the sheaf complex $rat(K)$. Similarly we use the notation $H_\chi^\bullet(X, K) = H^\bullet(X, rat(K)_\chi)$ for twists by characters $\chi$. Notice however that we do not assume that the character twisting functor lifts from the derived category $\Dbc(X, \Lambda)$ to the category~$\bfD$. Depending on the context we require some of the following three sets of axioms.

\medskip

\begin{itemize}
\item[(D1)] {\it Degree shifts}. We have an autoequivalence $K\mapsto K[1]$ on $\bfD$ which induces on $\Dbc(X, \Lambda)$ the usual degree shifting functor.
\medskip

\item[] {\it Perverse truncations}. We have endofunctors $\ptau_{\leq 0}, \ptau_{\geq 0}: \bfD \rightarrow \bfD$ and natural transformations $\ptau_{\leq 0} \to \id_\bfD \to \ptau_{\geq 0}$ which induce on~$\Dbc(X, \Lambda)$ the truncations for the perverse $t$-structure. 
\medskip

\item[] {\it Exactness}. The perverse cohomology functor  
$\pH^0 = \ptau_{\leq 0} \circ \ptau_{\geq 0}$
has as its essential image a full {\em abelian} subcategory $\bfP \subset \bfD$, and the given functor $rat: \bfP \rightarrow \Perv(X,\Lambda)$ is an exact functor from this abelian category to the abelian category of perverse sheaves.\medskip

\item[(D2)] {\it Semisimplicity}. For all objects $K\in \bfD$ there exists a (non-canonical) isomorphism 
\[ \quad K \;\;\cong \;\; \bigoplus_{n\in \bbZ} \; \pH^{n}(K)[-n]
\quad 
\textnormal{where} 
\quad
\pH^n(K) \;=\; \pH^0(K[n]). \] 
Furthermore, in axiom (D1) the abelian category $\bfP$ is semisimple. \medskip

\item[(D3)] {\it Hard Lefschetz}. In $\bfD$ there exists an invertible object ${\bf 1}(1)$ whose image in $\Perv(X,\Lambda)$ under $rat$ is the Tate twist of $\bf 1$. For all $K, L$ in~$\bfP$ and all $n\in \bbN$ we have functorial Lefschetz isomorphisms
$$  \pH^{-n}(K*L) \; \cong \; \pH^n(K*L)(n) ,$$
where the Tate twist $(n)$ denotes the $n$-fold convolution with $\one(1)$. \medskip
\end{itemize}

\noindent
Note that we do not assume $\bfD$ to be triangulated, indeed we will later deal with the following non-triangulated categories.  

\medskip

\begin{ex} \label{ex:pervgeom} 
The axioms (D1) -- (D3) are satisfied if $\bfD \subseteq \Dbc(X, \Lambda)$ is the full subcategory of all direct sums of degree shifts of semisimple perverse sheaves which in case $char(k)>0$ are defined over some finite field. 
\end{ex}

\medskip

Indeed, for $k=\bbC$ this holds by Kashiwara's conjecture, which has been reduced by Drinfeld~\cite{Dr} to a conjecture of de Jong that was proven some years later in~\cite{BK} and~\cite{Gai}. Alternatively, for $k=\bbC$ one can use the theory of polarizable twistor modules~\cite{Sab},~\cite{Moc}. In the case where $char(k)>0$ one can instead invoke the mixedness results of \cite{Laf}. Note that in the above example we could also replace the category $\bfD$ by the full subcategory  of objects of geometric origin in the sense of \cite[sect.~6.2.4]{BBD}. 

\medskip

\begin{ex} \label{ex:mhm}
The axioms (D1) -- (D3) hold for $k=\bbC$ and $\Lambda = \bbQ$ if $\bfD$ is taken to be the full subcategory of $D^b(\MHM(X))$ consisting of all direct sums of degree shifts of semisimple Hodge modules. 
\end{ex}

\medskip

For the proof of theorem~\ref{thm:P} we will consider a full subcategory~$\bfN$ of $\bfD$ consisting of objects that are negligible for our purposes. Since we want to proceed as far as possible over a base field of arbitrary characteristic, we formulate the required properties in the following axiomatic way.

\medskip

\begin{itemize}
 \item[(N1)] {\it Stability}. We have $\bfN * \bfD \subseteq \bfN$, and $\bfN$ is stable under taking direct sums, retracts, degree shifts, perverse truncations and adjoint duals.\medskip

 \item[(N2)] {\it Twisting}. Every object $K\in \bfN$ has the property $H^\bullet_\chi(X, K)=0$ for most characters $\chi$ of the fundamental group. \medskip
 
 \item[(N3)] {\it Acyclicity}. The category $\bfN$ contains all $K\in \bfD$ which are acyclic in the sense that $H^\bullet(X, K)=0$. \medskip
 
 \item[(N4)] {\it Euler characteristics}. The category $\bfN$ contains all simple objects of~$\bfP$ whose Euler characteristic vanishes. \medskip
\end{itemize}

\noindent
The meaning of these axioms will become clear later on. For the time being we content ourselves with the following

\medskip \begin{rem} Let $\Pi$ be a set of characters of $\pi_1(X, 0)$, and $\bfN\subseteq \bfD$ the full subcategory of all $K\in \bfD$ such that $rat(K)$ is a direct sum of degree shifts of local systems~$L_\chi$ with $\chi \in \Pi$. Then axioms (N1) and~(N2) hold. \end{rem} \medskip

{\em Proof.} For any $M\in \Dbc(X, \Lambda)$ we have $L_\chi * M = L_\chi \otimes_\Lambda H^\bullet(X, M_{\chi^{-1}})$ by~\cite[p.~20]{WeBN}, which in particular implies the stability property $\bfN*\bfD \subseteq \bfN$ so that axiom~(N1) holds. For (N2) use that $H^\bullet(X, L_\chi)=0$ if and only if the character $\chi$ is non-trivial.\qed

\section{The André-Kahn quotient} \label{sec:andre-kahn}

For our Tannakian arguments we want to work in rigid symmetric monoidal categories which are {\em semisimple abelian}. To construct such a category $\bfDbar$ which is as close as possible to the category $\bfD$, we use a general method of Andr\'e and Kahn~\cite{AK1} as explained below. In this section we always assume that the first two axioms (D1) and (D2) of section~\ref{sec:setting} hold.

\medskip

By rigidity, any endomorphism $f$ of an object $K$ in $\bfD$ has an adjoint morphism $f^\sharp: \one \to K*K^\vee$. The trace $tr(f)\in \End_\bfD(\one) = \Lambda$ is defined as the composite
$
 tr(f) = \mathrm{ev}_K \circ S_{K, K^\vee} \circ f^\sharp
$
where $S_{K, K^\vee}: K*K^\vee \to K^\vee*K$ denotes the symmetry constraint and where $\mathrm{ev}_K: K^\vee * K \to \one$ is the evaluation. As in section~7.1 of loc.~cit.~we consider the Andr\'e-Kahn radical~$\calN$ of $\bfD$, i.e.~the ideal which is defined for objects $K, L$ of $\bfD$ by
\[
 \calN(K, L) \;=\; \{ f\in \Hom_\bfD(K, L) \mid \forall g\in \Hom_\bfD(L, K): \; tr(g\circ f) = 0 \}.
\]
By definition, the quotient category
\[ \bfDbar = \bfD / \calN \]
has the same objects as $\bfD$, but the morphisms between two objects $K, L$ are defined by 
\[ \Hom_{\bfDbar}(K, L) \;=\; \Hom_{\bfD}(K, L)/\calN(K, L). \]
We have a natural quotient functor $q : \bfD \longrightarrow \bfD$ that is given by the identity on objects and by the quotient map on morphisms, and in what follows we denote by $\bfPbar$ the essential image of $\bfP$ under this quotient functor. Ultimately we want to construct a semisimple abelian category; as a first step towards this goal we have

\medskip \begin{lem} \label{lem:psab}
The quotient functor $q: \bfD \to \bfDbar$ preserves direct sums, and the category $\bfPbar$ is pseudo-abelian in the sense that every idempotent morphism in it splits as the projection onto a direct summand.
\end{lem} \medskip

{\em Proof.} The functor $q$ preserves direct sums since it is $\Lambda$-linear. To see that idempotents in $\bfPbar$ split, let $P$ be an object of $\bfP$. Since $\bfP$ is an abelian category, it suffices to show that every idempotent in
 \[\End_\bfPbar(P) \;=\; \End_\bfP(P)/\calN(P, P)\]
lifts to an idempotent in $\End_\bfP(P)$.
Since $\bfP$ is semisimple by axiom~(D2), we can assume $P=Q^{\oplus r}$ for some simple object~$Q$ of $\bfP$ and $r\in \bbN$. Then $\End_\bfP(P)$ is the ring of $r\times r$ matrices over the skew field $\End_\bfP(Q)$. Since matrix rings over skew fields do not have proper two-sided ideals, it follows that either $\calN(P, P)=0$ or $\calN(P,P)=\End_\bfP(P)$. In both cases the lifting of idempotents is obvious. \qed

\medskip

\medskip \begin{prop} \label{prop:abelian}
The quotient category $\bfDbar$ is a $\Lambda$-linear semisimple abelian rigid symmetric monoidal category. 
\end{prop} \medskip

{\em Proof.} By lemma~7.1.1 in loc.~cit.~the Andr\'e-Kahn radical $\calN$ is a monoidal ideal, so it follows from sorite 6.1.4 of loc.~cit.~that the quotient category $\bfDbar$ is again a \mbox{$\Lambda$-linear} rigid symmetric monoidal category with $End_{\bfDbar}({\bf 1})=\Lambda$. We claim that
\begin{equation} \label{eq:HomVanishing} 
\Hom_\bfDbar(P[m], Q[n]) \;=\; 0 \quad \textnormal{for all $P, Q$ in $\bfP$ and $m\neq n$}.
\end{equation}
Indeed, for $m> n$ we even have $\Hom_\bfD(P[m], Q[n]) = 0$ since under the faithful functor $rat$ this $\Hom$-group injects into 
\[
 \Hom_{\Dbc(X, \Lambda)} (rat(P)[m], rat(Q)[n]) \;=\; \Ext_{\Perv(X, \Lambda)}^{\, n-m}(rat(P), rat(Q))
\] 
which vanishes for $m>n$ (for the above interpretation as an $\Ext$-group recall that $\Dbc(X, \Lambda)$ is the derived category of $\Perv(X, \Lambda)$). For $m<n$ similarly $\Hom_\bfD(Q[n], P[m])=0$, and in that case the definition of $N_\bfD$ trivially implies that $\Hom_\bfD(P[m], Q[n]) = N_\bfD(P[m],Q[n])$. This is mapped to zero under the quotient functor $\bfD \to \bfDbar$, hence our claim~\eqref{eq:HomVanishing} follows. 

\medskip

Now by the semisimplicity axiom~(D2) every object $K$ of $\bfDbar$ can be written as $K=\bigoplus_{n\in \bbZ} K_n[n]$ with certain $K_n$ in~$\bfPbar$. The vanishing property in \eqref{eq:HomVanishing} then implies
\smallskip
\begin{equation} \label{eq:End} 
\End_\bfDbar(K) \;=\; \bigoplus_{n\in \bbZ} \, \End_\bfDbar (K_n[n]) \;=\; \bigoplus_{n\in \bbZ} \, \End_\bfPbar (K_n).
\end{equation}
In particular, every idempotent endomorphism of $K$ in the category $\bfDbar$ is a direct sum of idempotent endomorphisms of the summands $K_n[n]$, and by lemma~\ref{lem:psab} it follows that $\bfDbar$ is pseudo-abelian. Hence to show that $\bfDbar$ is a semisimple abelian category, it will suffice by \cite[A.2.10]{AK1} to show that it is a semisimple $\Lambda$-linear category in the sense of section 2.1.1 in loc.~cit. For this we use the following general result~\cite[th.~1]{AK2}:

\medskip

{
Let $F$ be a field and $\bfA$ an $F$-linear rigid symmetric monoidal category with $\End_\bfA(\one)=F$. Suppose there exists an $F$-linear tensor functor ACU from~$\bfA$ to an abelian $F$-linear rigid symmetric monoidal category $\bfV$ such that $\dim_\Lambda ( \Hom_\bfV (V_1, V_2)) < \infty$ for all $V_1, V_2 \in \bfV$. Then the quotient of~$\bfA$ by its Andr\'e-Kahn radical $N_\bfA$ is a semisimple $F$-linear category, and $N_\bfA$ is the unique monoidal ideal of $\bfA$ with this property.}

\medskip

In our case this applies for $F=\Lambda$, $\bfA=\bfD$ and for the functor $H^\bullet(X, -)$ from $\bfD$ to the abelian category $\bfV$ of super vector spaces over $\Lambda$. \qed

\smallskip

\medskip \begin{cor} \label{cor:exact}
The functors $\bfP \to \bfPbar$ and $\bfPbar \hookrightarrow \bfDbar$ are exact functors between semisimple abelian categories. The image of a simple object $P\in \bfP$ inside $\bfPbar$ is either simple or isomorphic to zero, and if $\Lambda$ is algebraically closed, then the latter case occurs if and only if $\chi(P)=0$.
\end{cor} \medskip

{\em Proof.} By proposition~\ref{prop:abelian}, $\bfDbar$ is a semisimple abelian category, and it also follows from the proof of the proposition that $\bfPbar$ is a semisimple abelian subcategory of $\bfDbar$. Since the considered functors are additive, they are exact by semisimplicity. If~$P$ is a simple object of $\bfP$, then $\End_\bfP(P)$ is a skew field, hence $\End_\bfPbar(P)$ is a skew field or zero, and $P$ is simple or zero in $\bfPbar$. Over an algebraically closed field $\Lambda$ there exist no skew fields other than $\Lambda$ itself, hence in this case we have $\End_\bfP(P)=\Lambda$, and it follows that $\id_P \in N_\bfD(P, P)$ iff $tr(\id_P)=0$, which is the case iff $\chi(P)=0$.
\qed

\smallskip

\medskip \begin{cor} \label{cor:iso_to_zero} Let $\bfN\subseteq \bfD$ be the full subcategory of all objects which become isomorphic to zero in the quotient category $\bfDbar$. If $\Lambda$ is algebraically closed, then $\bfN$ satisfies the stability axiom (N1), the acyclicity axiom (N3) and the Euler axiom (N4), and an object $K\in \bfD$ lies in the subcategory $\bfN$ iff all simple constituents of all $\pH^n(K)$ have Euler characteristic zero. 
\end{cor} \medskip

{\em Proof.} Property (N1) is clear, (N3) follows from (N4), and the latter is immediate from corollary~\ref{cor:exact} in view of the semisimplicity axiom (D2). \qed

\section{Super Tannakian categories} \label{sec:tannakian}

Using a criterion of Deligne, we now show that the semisimple abelian rigid symmetric monoidal category $\bfDbar$ from the previous section is almost Tannakian: It is an inductive limit of finitely generated super Tannakian categories, a notion that we will recall below and in the appendix. For $k=\bbC$ we will see in corollary~\ref{group} that $\bfDbar$ is an inductive limit of finitely generated ordinary Tannakian categories, a fact closely related to theorem~\ref{thm:P}. 

\medskip

In this section we always assume that $\Lambda$ is algebraically closed and that the first two axioms (D1) and (D2) from section~\ref{sec:setting} are satisfied. By semisimplicity the functor $*: \bfDbar\times \bfDbar \to \bfDbar$ is exact in each variable, and $\End_\bfDbar(\one)=\Lambda$ (this is inherited from $\bfD$ and can be checked via the faithful functor $rat$). Hence~$\bfDbar$ is a {\em cat\'egorie $\Lambda$-tensorielle} in the sense of~\cite[sect.~0.1]{DelCT}. 

\medskip

Recall that a full subcategory of $\bfDbar$ is said to be finitely tensor generated, if it is the category of all subquotients of convolution powers of~$C\oplus C^\vee$ for some fixed object $C$. 
The next theorem will show that any such category is super Tannakian in the following sense.

\medskip

The framework of algebraic geometry can be generalized to super algebraic geometry by replacing the category of commutative rings with the one of $\bbZ/2\bbZ$-graded super commutative rings. In particular one has the notions of algebraic and reductive super groups over $\Lambda$ and their super representations, as we recall in the appendix in section~\ref{sec:supergroups} below. For an algebraic super group~$G$ over~$\Lambda$ and a point $\epsilon \in G(\Lambda)$ with $\epsilon^2 = 1$ such that $\mathrm{int}(\epsilon)$ is the parity automorphism of $G$, we denote by $\Rep_\Lambda(G, \epsilon)$ the category of super representations $V=V_+\oplus V_-$ of $G$ over $\Lambda$ for which $\epsilon$ acts by $\pm 1$ on $V_\pm$.
Such categories will be called super Tannakian with Tannaka super group $G$.

\medskip \begin{thm} \label{thm:deligne}
Every finitely generated full tensor  subcategory $\bfT$ of $\bfDbar$ is super Tannakian with a reductive Tannaka super group $G=G(\bfT)$. 
\end{thm} \medskip

{\em Proof.} Since $\bfDbar$ is a cat\'egorie $\Lambda$-tensorielle, for the first claim it suffices by \cite[th.~0.6]{DelCT} to see that for any object $C\in \bfDbar$ the number of constituents of $C^{*n}$ is at most $N^n$ for some constant $N=N(C)$ and all $n\in \bbN$. For this one can take $N(C)=\sum_{i\in \bbZ} \dim_\Lambda (H^i(X, D))$ with any object $D\in \bfD$ that becomes isomorphic to $C$ in $\bfDbar$, see~\cite[top of p.~5]{WeT}.
Concerning reductivity, note that by~\cite{WeSS} a category $\Rep_\Lambda(G, \epsilon)$ is semisimple iff $G$ is reductive. 
\qed

\section{Perverse multiplier} \label{sec:multiplier}

We now introduce the notion of a perverse multiplier with respect to a given subcategory of negligible objects; this notion will play an important role in our proof of theorem~\ref{thm:P}. In this section $\Lambda$ need not be algebraically closed, but we still assume that axioms (D1) and (D2) of section~\ref{sec:setting} hold, and we consider a full subcategory $\bfN \subseteq \bfD$ with the stability properties (N1).

\medskip \begin{defn} \label{def:multiplier}
An object $K\in \bfD$ is called an $\bfN$-multiplier, if for all $r\in \bbN_0$ and all $n\neq 0$ every subquotient of~$\pH^n((K\oplus K^\vee)^{*r})$ lies in $\bfN$. We say that~$K$ is a zero type, if $H^n(X, K)=0$ holds for all $n\neq 0$.
\end{defn} \medskip

The relevance of these notions for the proof of theorem~\ref{thm:P} becomes clear from the following observation.

\medskip \begin{lem} \label{lem:mult-type}
For $P\in \bfP$ the following holds. 
\begin{itemize} 
 \item[\em (a)] If $\bfN$ satisfies the twisting axiom (N2) and if~$P$ is an $\bfN$-multiplier, then $H^\bullet_\chi(X, P)$ is concentrated in degree zero for most $\chi$. \smallskip 
 \item[\em (b)] If $\bfN$ satisfies the acyclicity axiom (N3), if the hard Lefschetz axiom~(D3) holds and if $P$ is a zero type, then $P$ is an $\bfN$-multiplier.
\end{itemize}
\end{lem} \medskip

{\em Proof.} {\em (a)} Put $g=\dim(X)$. The semisimplicity axiom (D2) shows that we then have
\vspace*{0.1em}
\[ P^{*(g+1)} \;=\; \bigoplus_{m\in \bbZ} \; P_m[m] \quad \textnormal{for suitable} \quad P_m \in \bfP \,. \]
By assumption $P$ is an $\bfN$-multiplier, hence $P_m\in \bfN$ for all~$m\neq 0$. Via the twisting axiom (N2) it follows that for most characters $\chi$ and all $n\in \bbZ$,
\[
 H^n_\chi(X, P^{*(g+1)}) \;=\; H^n_\chi(X, P_0).
\]
The right hand side vanishes for $|n|>g$, since $rat(P_0)_\chi$ is perverse. But for the left hand side we have
\[
 H^\bullet_\chi(X, P^{*(g+1)})  \;=\; (H^\bullet_\chi(X, P))^{\otimes (g+1)} 
\]
by proposition~\ref{prop:tensorfunctor} and since $H^\bullet(X, -)$ is a tensor functor by the K\"unneth theorem. So the above vanishing statement for~$|n|>g$ implies that $H^\bullet_\chi(X,P)$ is concentrated in degree zero. 

\medskip

{\em (b)} Put $Q=(P\oplus P^\vee)^{*r}$ for any $r\in \bbN$. Since hypercohomology is a tensor functor by the K\"unneth theorem, with $H^\bullet(X, P)$ also $H^\bullet(X, Q)$ is concentrated in degree zero. Using the hard Lefschetz axiom (D3), one then deduces that for all $n\neq 0$ one has $H^\bullet(X,\pH^n(Q))=0$ so that by (N3) the subcategory $\bfN$ contains $\pH^n(Q)$. Since this holds for arbitrary $r\in \bbN$, it follows that indeed $P$ is an $\bfN$-multiplier.  \qed

\medskip

In view of part {\em (a)} of the lemma, to prove theorem~\ref{thm:P} we want to show that for a suitable subcategory $\bfN$ every object of $\bfP$ is an $\bfN$-multiplier. For this we will argue by contradiction, using the following

\medskip \begin{lem} \label{lem:skyscraper}
Suppose that $\bfN$ satisfies the stability axiom (N1) and the Euler axiom (N4), that $\bfD$ satisfies all axioms (D1) -- (D3) and that $P\in \bfP$ is not an $\bfN$-multiplier. Then for some $r\in \bbN$ the convolution power 
\[ (P * P^\vee)^{*r} \;=\; (P*P^\vee)*\cdots * (P*P^\vee) \] 
admits a direct summand of the form ${\bf 1}[2i](i)$ with an integer~$i\neq 0$.
\end{lem} \medskip

{\em Proof.} If $P$ is not an $\bfN$-multiplier, we can find integers $a, b\in \bbN$ such that~$P^{*a} * (P^\vee)^{*b}$ admits a direct summand $Q[i]$ for some $i\neq 0$ and some simple object $Q\in \bfP$ which is not in $\bfN$. By the hard Lefschetz axiom (D3) then $Q[-i](-i)$ is a direct
summand of $P^{*a} * (P^\vee)^{*b}$ as well. It then follows that also the dual $Q^\vee[i](i)$ is a direct summand of $P^{*b}*(P^\vee)^{*a}$. Altogether then the convolution product $Q[i]*Q^\vee [i](i) = Q*Q^\vee [2i](i)$ will be a direct summand of $(P * P^\vee)^{*r}$ for the exponent $r=a+b$.

\medskip

It remains to show that $\one$ is a direct summand of $Q*Q^\vee$. For this note that the trace map $tr(Q): \one \longrightarrow Q*Q^\vee \cong Q^\vee * Q \longrightarrow \one$ is non-zero, since we have $\chi(Q)\neq 0$ by axiom (N4). Now $tr(Q)$ factors over $\pH^0(Q*Q^\vee)$, indeed $\Hom_\bfD(\bfP, \ptau_{>0} \bfP)=\Hom_\bfD(\ptau_{<0} \bfP, \bfP)=0$. So $tr(Q)$ exhibits $\one$ as a retract of $\pH^0(Q*Q^\vee)$ in the abelian category $\bfP$, and we are done. 
\qed

\medskip

\section{Proof of the vanishing theorem} \label{sec:proof}

The main idea of our proof of theorem~\ref{thm:P} is to control the non-perversity of convolution products in terms of central characters of the Tannaka group from theorem~\ref{thm:deligne}. By d\'evissage we can  restrict ourselves to semisimple perverse sheaves as in example~\ref{ex:pervgeom}. So suppose that $\Lambda = \bbC$ or $\Lambda = \Qbar_l$ and that~$\bfD$ satisfies all axioms (D1) -- (D3) of section~\ref{sec:setting}. Consider the semisimple abelian rigid symmetric monoidal quotient category $\bfDbar$ from section~\ref{sec:andre-kahn}.

\medskip

For the full subcategory $\bfN\subseteq \bfD$ of all objects that become isomorphic to zero in~$\bfDbar$, the axioms (N1), (N3) and (N4) hold by corollary~\ref{cor:iso_to_zero}. We expect that in the setting of example~\ref{ex:pervgeom} also axiom (N2) always holds. However, at present we can show this only for $k=\bbC$ via the theory of $\calD$-modules, which we will do in corollary~\ref{cor:acyclictwist} below. In any case, once we have (N2), we can apply part (a) of lemma~\ref{lem:mult-type} to deduce the vanishing theorem~\ref{thm:P} from the axioms (N1) and (N4) via the next

\medskip \begin{thm} \label{thm:mult}
Let $\bfN\subseteq \bfD$ be a full subcategory satisfying the axioms (N1) and (N4). Then every object $P\in \bfP$ is an $\bfN$-multiplier.
\end{thm} \medskip

{\em Proof.} Suppose that $P\in \bfP$ is simple and not an $\bfN$-multiplier.  Then for some integer $r\in \bbN$ the convolution $(P * P^\vee)^{*r}$ contains by lemma~\ref{lem:skyscraper} a direct summand $L={\bf 1}[2i](i)$ with $i\neq 0$. Hence the full rigid symmetric monoidal subcategory $\bfDbar_1$ generated inside $\bfDbar$ by $P$ contains the full rigid symmetric monoidal subcategory $\bfDbar_0$ generated by the invertible object $L$. 

\medskip

Theorem~\ref{thm:deligne} shows that for certain reductive super groups~$G_i$ over $\Lambda$ we have tensor equivalences $\omega_i: \bfDbar_i \stackrel{\sim}{\longrightarrow} \Rep_\Lambda(G_i, \epsilon_i)$ for $i\in \{0,1\}$, and by the Tannakian formalism the inclusion $\bfDbar_0\subseteq \bfDbar_1$ defines an epimorphism of reductive super groups 
$$h: \; G_1\twoheadrightarrow G_0.$$ 
The category $\bfDbar_0$ consists of all direct sums of skyscraper sheaves $L^{*n}$ with integers $n\in \bbZ$. Since
$L^{*n}\cong {\bf 1}[2ni](ni)$ and $i\neq 0$, equation~\eqref{eq:End} in the proof of proposition~\ref{prop:abelian} implies that one has $L^{*n} \cong {\bf 1}$ in $\bfDbar$ only if $n=0$. Taking into account that the symmetry constraint  
$L*L \longrightarrow L*L$ is the identity in $\bfD$, the tensor equivalence $\omega_0$ between $\bfDbar_0$ and $\Rep_\Lambda(\bbG_m, 1)$ is realized explicitly, with the multiplicative Tannaka group $G_0=\bbG_m$ and $\epsilon_0=1$,  via
$$  \quad L^{*n} \ \mapsto \  (\textnormal{the character $z\mapsto z^{n}$ of $\bbG_m$}).
$$
In particular, the representation $W_0=\omega_0(L)$ is non-trivial.

\medskip

But proposition~\ref{prop:toruslift} in the appendix applies to the torus $T_0 = G_0 = \bbG_m$, so there exists a central torus $T_1\cong \bbG_m$ in $G_1$ such that $h: G_1\to G_0$ restricts to an isogeny $T_1\to T_0$. By Schur's lemma the central torus $T_1$ acts via some character on the irreducible super representation $W_1=\omega_1(P)$; so $T_1$ acts trivially on~$W_1\otimes W_1^\vee = \omega_1(P*P^\vee)$. Then $T_1$, hence also $T_0$,  acts trivially also on the direct summand $W_0\subseteq (W_1\otimes W_1^\vee)^{\otimes r}$ -- a contradiction. \qed

\medskip \begin{cor}\label{group}
In the case of the base field $k=\bbC$, the super group $G(\bfT)$ in theorem~\ref{thm:deligne} is a classical reductive algebraic group over $\Lambda$.
\end{cor} \medskip

{\it Proof}. Corollary \ref{cor:iso_to_zero} and theorem~\ref{thm:mult} show that the category $\bfPbar$ is preserved
under convolution. Using this one easily reduces our claim to the special case where $\bfT \subset \bfPbar$. The assertion then follows from \cite[th.~7.1]{DeTann} since for $k=\bbC$ we will see in section~\ref{sec:euler} that $\chi(P)\geq 0$ for all $P\in \bfPbar$. 
\qed

\section{Euler characteristics} \label{sec:euler}

In view of corollary~\ref{cor:iso_to_zero}, to control how much information is lost in the passage from $\bfD$ to $\bfDbar$ we must determine all perverse sheaves on $X$ with Euler characteristic zero. This will complete the proof of theorem~\ref{thm:P}, since it will imply that the category $\bfN$ from section~\ref{sec:proof} satisfies axiom (N2). In this section we always work over $k=\bbC$. Then by~\cite[cor.~1.4]{FK} every perverse sheaf $P$ has Euler characteristic $\chi(P)\geq 0$, and we have

\medskip 

\begin{prop} \label{prop:chi}
Let $P$ be a simple perverse sheaf on $X$.
\begin{itemize}
\item[\em (a)] One has $\chi(P)=0$ iff there exists a positive-dimensional abelian subvariety $A\hookrightarrow X$ with quotient $q: X\twoheadrightarrow B=X/A$ such that
\[
 P \;\cong \; L_\varphi \otimes q^*(Q)[\dim(A)]
\]
for some $Q\in \Perv(B, \bbC)$ and some character $\varphi$ of $\pi_1(X, 0)$.  \medskip

\item[\em (b)] One has $\chi(P)=1$ iff $P$ is a skyscraper sheaf on $X$ of rank one.
\end{itemize}
\end{prop} \medskip

{\em Proof.} View $P$ as a $\calD_X$-module via the Riemann-Hilbert correspondence. For $Z\subseteq X$ closed and irreducible, let $\Lambda_Z \subseteq T^*X$ be the closure in $T^*X$ of the conormal bundle in $X$ to the smooth locus of~$Z$. As in loc.~cit.~we write the characteristic cycle of $P$ as a finite formal sum
\smallskip
\[
 \CC (P) \;=\; \sum_{Z\subseteq X} n_Z\cdot \Lambda_Z
 \quad \textnormal{with} \quad n_Z\in \bbN_0,
\]
where $Z$ runs through all closed irreducible subsets of $X$.
From $\CC(P)$ the support of the perverse sheaf $P$ can be recovered via
$\Supp \, P \,  =  \bigcup_{n_Z\neq 0} Z$.
Furthermore, by the microlocal index formula~\cite[th.~9.1]{Gi}, 
\smallskip
\[
 \chi(P) \;=\; \sum_{Z\subseteq X} n_Z\cdot d_Z \quad \textnormal{with} \quad d_Z = [\Lambda_X]\cdot [\Lambda_Z] \in \bbZ.
\]
The intersection numbers $d_Z$ are well-defined even though $\Lambda_Z$ is not proper for $Z\neq X$; see loc.~cit.~for details. Now if $X$ is a simple abelian variety, then lemma~\ref{lem:numbers} below implies the claim {\em (a)}, and if we additionally assume $\dim(X)>1$, also {\em (b)} follows in view of lemma~\ref{lem:cc} below. The non-simple case can be reduced to the simple case, see~\cite{WeGauss}.
\qed 

\medskip

The reduction step to the case of simple abelian varieties in~\cite{WeGauss} works for ground fields $k$ of characteristic $p>0$ as well, but for (simple) abelian varieties defined over a finite field the above argument has to be replaced by a kind of Iwasawa-theoretic deformation argument~\cite{WeVanishing}. For $k=\bbC$, Christian Schnell has given in~\cite[th.~7.6]{Schn} a different proof of proposition~\ref{prop:chi}(a) using the Fourier-Mukai transform for $\calD_X$-modules.

\medskip \begin{cor} \label{cor:acyclictwist}
The Euler characteristic of a simple perverse sheaf $P$ on~$X$ vanishes iff
$H^\bullet(X, P_\chi)=0$
for most characters $\chi$.
\end{cor} \medskip

{\em Proof.} ``$\Leftarrow$'' holds by corollary~\ref{cor:chi}. For ``$\Rightarrow$'' take a positive-dimensional abelian subvariety $A\hookrightarrow X$ with quotient $q: X\twoheadrightarrow B=X/A$ and a character~$\varphi$ such that $P \cong L_\varphi \otimes q^*(Q)[\dim(A)]$ for some perverse sheaf $Q$ on $B$ as in proposition~\ref{prop:chi}a). We can assume that the Euler characteristic of $Q$ is not zero. Then we claim that
\[ H^\bullet(X, P_\chi) = H^\bullet(B, Rq_*(P_\chi)) = H^\bullet(Rq_*(L_{\varphi \chi}) \otimes Q[\dim(A)]) \] 
vanishes iff the restriction of the local system $L_{\varphi \chi}$ to $A=\ker(q)$ is not trivial. Indeed, if this restriction is non-trivial, then $Rq_*(L_{\varphi \chi})=0$ and hence also $H^\bullet(X, P_\chi)=0$. But if this restriction is trivial, then $L_{\varphi \chi} = q^*(L_\psi)$ for some character $\psi$, and then $H^\bullet(X, P_\chi) = H^\bullet(A, \bbC) \otimes H^\bullet(B, Q_\psi)[\dim(A)]$ is non-zero since the Euler characteristic of $Q_\psi$ is not zero.
\qed

\medskip \begin{lem} \label{lem:numbers}
One has $d_Z\geq 0$ for all $Z$. Furthermore, $d_Z = 1$ iff $Z$ is reduced to a single point. If $X$ is simple, then $d_Z = 0$ iff $Z=X$.
\end{lem} \medskip

{\em Proof.} The cotangent bundle $T^*X = X\times \bbC^g$ is trivial of rank $g=\dim(X)$, and projecting from $\Lambda_Z\subseteq T^*X$ onto the second factor $\bbC^g$ induces the Gau\ss\ mapping $p: \Lambda_Z \to \bbC^g$. By~\cite[prop.~2.2]{FK} the intersection number $d_Z$ is the generic degree of $p$. In particular $d_Z\geq 0$.

\medskip

If $d_Z=1$, then $\Lambda_Z$ is birational to $\bbC^g$, so by~\cite[cor.~3.9]{Mi} there does not exist any non-constant map from $\Lambda_Z$ to an abelian variety. So the image $Z$ of the composite map $\Lambda_Z\subseteq T^*X \to X$ is a single point. 

\medskip

If $d_Z=0$, then $p$ is not surjective, so $\dim(p(\Lambda_Z)) < g$. Then for some cotangential vector $\omega \in p(\Lambda_Z)$ the fibre $p^{-1}(\omega)$ is positive-dimensional. If $Z\neq X$, we can assume $\omega \neq 0$. Let $Y\subseteq X$ be the image of $p^{-1}(\omega)\subseteq T^*X$ under the map $T^*X \to X$. Then $\dim(Y)>0$, and up to a translation we can assume $0\in Y$. By construction $\omega$ is normal to $Y$ in  every smooth point of~$Y$, so the preimage of $Y$ under the universal covering $\bbC^g\to X=\bbC^g/\Lambda$ lies in the hyperplane of $\bbC^g$ orthogonal to $\omega$. Thus the abelian subvariety of~$X$ generated by $Y$ is strictly contained in $X$ but non-zero, contradicting the assumption that $X$ is simple.
\qed

\medskip

\medskip \begin{lem} \label{lem:cc}
Let $P$ be a simple perverse sheaf on $X$. If there is a closed subset $Y\subset X$ with $\dim(Y)\leq g-2$ such that %
\smallskip
\[ \CC (P) \;=\; n_X \Lambda_X + \sum_{Z\subseteq Y} n_Z \Lambda_Z \quad \textnormal{\em and} \quad n_X > 0, \]
then $P = L_\chi[g]$ for the local system $L_\chi$ on $X$ attached to some character $\chi$.
\end{lem} \medskip

{\em Proof.} Consider the embedding $j: U=X\setminus Y\hookrightarrow X$. Open embeddings are non-characteristic for any $\calD_X$-module, so theorem 2.4.6 and remark 2.4.8 in~\cite{HTT} show
$
 CC(j^*(P)) = CC(P)\cap T^*U = n_X\cdot \Lambda_U.
$
By prop.~2.2.5 in loc.~cit.~then $j^*(P)=L_U[g]$ for some local system $L_U$ on~$U$. Since $X$ is smooth, by the purity of the branch locus the assumption on $\dim(Y)$ implies $L_U=j^*(L)$ for some local system $L$ on $X$. By simplicity of~$P=j_{!*}(j^*(P))$ then $L$ has rank one, and $P=L[g]$. 
\qed

\section{The spectrum of a perverse sheaf} \label{sec:spectrum}

Let $X$ be a complex abelian variety of dimension $g$. Then $\pi_1(X,0) \cong \bbZ^{2g}$ and the group $\Pi(X)$ of characters $\chi: \pi_1(X,0) \to \bbC^*$ is a complex algebraic torus of rank $2g$. For any semisimple perverse sheaf $P$ on $X$, we explain in this section how to determine the set of all $\chi \in \Pi(X)$ for which theorem~\ref{thm:P} fails, and in particular we show that this set is a finite union of translates of proper algebraic subtori of $\Pi(X)$. We also consider the corresponding question in the relative setting of theorem~\ref{thm:relative}.

\medskip

Note that $\Pi$ is a contravariant functor: Any homomorphism $h: X\to B$ of abelian varieties induces a homomorphism $\pi_1(h): \pi_1(X, 0)\to \pi_1(B, 0)$ and hence a homomorphism $\Pi(h): \Pi(B) \to \Pi(X)$ of algebraic tori. For a perverse sheaf $P$ on $X$ we define the spectrum ${\calS}(P) \subset \Pi(X)$ to be the set of all $\chi \in \Pi(X)$ such that
\[ H^i(X,P_\chi) \ \neq \ 0 \quad \textnormal{for some} \quad i\neq 0. \]
More generally, for a semisimple complex $K = \bigoplus_{n\in \bbZ} \pH^{-n}(K)[n]$ on $X$ we define
\[
 \calS(K) \;=\; \bigcup_{n\in \bbZ} \, \calS(\pH^{-n}(K)).
\]
It follows from the definitions that ${\calS}(K_\chi)=\chi^{-1}\cdot {\calS}(K)$ for all $\chi \in \Pi(X)$ and
that for all semisimple $K_1, K_2$ we have
\[  
 {\calS}(K_1 * K_2) \; \subseteq \;
 {\calS}(K_1) \cup {\calS}(K_2) \; = \; 
 {\calS}(K_1\oplus K_2).
\]
In particular, the last equality reduces the computation of the spectrum of semisimple sheaf complexes to the case of simple perverse sheaves. Note that $\calS(P)$ may be empty; for example, this is the case if $P$ is a skyscraper sheaf or if~$P=i_*E[1]$ where $i: C\hookrightarrow X$ is the embedding of a smooth curve in $X$ and where $E$ is an irreducible local system on $C$ of rank at least two.

\medskip

\medskip \begin{rem} \label{rem:Pi} The functor $\Pi$ has the following properties.
\begin{itemize} 
 \item[\em (a)] Let $g:X\to B$ be an isogeny with kernel $F$. Then we have an exact sequence
\[
\quad
\xymatrix{ 0 \ar[r] & \Hom(F, \bbC^*) \ar[r] & \Pi(B) \ar[r]^{\Pi(g)} & \Pi(X) \ar[r] & 0. 
}
\]
For any perverse sheaf $P$ on $X$ the direct image $g_*(P)$ is a perverse sheaf on $B$, and $\Pi(g)$ induces a surjection 
\[ \calS(g_*(P)) \twoheadrightarrow \calS(P). \] 
\item[\em (b)] Let $i: A\hookrightarrow X$ be the inclusion of an abelian subvariety with quotient morphism $q: X\to B=X/A$. Then we have an exact sequence 
\[ 
\xymatrix{ 0 \ar[r] & \Pi(B) \ar[r]^{\Pi(q)} & \Pi(X) \ar[r]^{\Pi(i)} & \Pi(A) \ar[r] & 0. 
}
\]
In this situation we denote by $K(A)\subseteq \Pi(X)$ the image of $\Pi(q)$.
\end{itemize}
\end{rem} \medskip

{\em Proof.} The exactness of the considered sequences can be seen from the description of a complex abelian variety as the quotient of a complex vector space modulo a lattice. For the surjectivity $\calS(g_*(P)) \twoheadrightarrow \calS(P)$
in part {\em (a)} use that $H^i(X,P_{\, \Pi(g)(\chi)})=H^i(B, g_*(P)_{\chi})$ and that $\Pi(g)$ is surjective. \qed

\medskip

In what follows, we denote by $E(X)$ the class of all semisimple perverse sheaves on $X$ with Euler characteristic zero. A perverse sheaf will be called clean, if it does not contain constituents from $E(X)$. For $x\in X(\bbC)$ we denote by $t_x: X\to X$ the translation morphism $y\mapsto x+y$, and for $K\in \Dbc(X, \bbC)$ we consider the stabilizer
\[
 \Stab(K) \;=\; \{ x\in X(\bbC) \mid t_x^* (K) \cong K\}.
\]
Its connected component $\Stab(K)^0 \subseteq \Stab(K)$ is an abelian subvariety of $X$.

\medskip \begin{lem} \label{cor:badlocus}
With notations as above, the following properties hold for the spectrum of semisimple perverse sheaves. 
\begin{enumerate}
\item[\em (a)] For $P\in E(X)$ we have $\calS(P) = \{ \chi \! \mid H^\bullet(X, P_\chi)\neq 0 \}$, and if $P$ is simple, there exists a character $\varphi$ such that
\[
 \calS(P) \;=\; \varphi^{-1} \cdot K(A) \quad \textnormal{\em for} \quad A = \Stab(P)^0
\]
where $K(A)\subset \Pi(X)$ is the proper subtorus from remark~\ref{rem:Pi}(b).

\smallskip

\item[\em (b)] For every semisimple $P\in \Perv(X, \bbC)$ there exist non-zero abelian subvarieties $A_i \subseteq X$ and characters $\chi_i\in \Pi(X)$ for $1\leq i\leq n$ with
\[
 \calS(P) \;=\; \bigcup_{i=1}^n \; \chi_i \cdot K(A_i).
\]
\item[\em (c)] If in part (b) the perverse sheaf $P$ is of geometric origin in the sense of~\cite[6.2.4]{BBD}, then the~$\chi_i$ can be chosen to be torsion characters.
\end{enumerate}
\end{lem} \medskip
 
{\em Proof.} {\em (a)} The first statement holds by proposition~\ref{prop:chi}(a), and the second one follows easily from the proof of corollary~\ref{cor:acyclictwist}.

\medskip

{\em (b)} By theorem~\ref{thm:mult} applied to the class $\bfN = \bfN_{Euler}$ of complexes with perverse cohomology sheaves in $E(X)$, we have
\[
 P^{*g} \;=\; Q \, \oplus \, \bigoplus_{\nu} \, N_\nu[\nu ]
\]
where $Q$ is a clean semisimple perverse sheaf and the $N_\nu$ are semisimple perverse sheaves in $E(X)$. Since twisting with a character is a tensor functor by proposition~\ref{prop:tensorfunctor}, it follows for any $\chi\in \Pi(X)$ that
\begin{equation} \label{eq:gthPower} 
 H^\bullet(X, P_\chi)^{\otimes g} \;=\; H^\bullet(X, Q_\chi) \, \oplus \, \bigoplus_{\nu} \, H^\bullet(X, (N_{\nu})_\chi)[\nu ].
\end{equation}
If $\chi \in \calS(P)$, then $H^\bullet(X, P_\chi)$ is not concentrated in degree zero, so~\eqref{eq:gthPower} is non-zero in some degree $d$ with $|d|\geq g$. But $Q_\chi$ is a clean perverse sheaf and as such it does not contain the constant perverse sheaf $\bbC[g]$ as a constituent, hence we have
\[
 H^d(X, Q_\chi) \;=\; 0 \quad \textnormal{for} \quad |d|\geq g
\]
by the adjunction properties in \cite[prop.~4.2.5]{BBD}. Thus $H^\bullet(X, (N_\nu)_\chi)\neq 0$ for some $\nu$ and hence $\chi \in \calS(N_{\nu})$ by part {\em (a)}. Conversely, if $\chi \in \calS(N_\nu)$ for some $\nu$, then $H^\bullet(X, (N_\nu)_\chi)$ is non-zero in at least two different cohomology degrees; then by~\eqref{eq:gthPower} the same holds for $H^\bullet(X, P_\chi)$, so $\chi \in \calS(P)$. Hence we have shown that
\[
 \calS(P) \;=\; \bigcup_{\nu} \; \calS(N_{\nu}),
\]
and our claim follows from the second statement in part {\em (a)}. 

\medskip

{\em (c)} First we claim that a local system $L_\chi$ is of geometric origin iff $\chi$ is a torsion character. For the non-trivial direction note that if $L_\chi$ is of geometric origin, then $X$ has a model~$X_A$ over a subring $A\subset \bbC$ of finite type over $\bbZ$ such that $L$ descends to a local system on $X_A$.
Take a closed point of $\mathit{Spec}(A)$ with finite residue field $\kappa$. Let $V\subset \bbC$ be a strictly Henselian ring with $A\subset V$ whose residue field is an algebraic closure $\kappabar$ of~$\kappa$. For $X_V = X_A \times_A V$ the inclusion of the special fibre $X_\kappabar$ induces an epimorphism $\pi_1(X_\kappabar) \twoheadrightarrow \pi_1(X_V)$ by the homotopy sequence \cite[exp.~X, 1.6]{SGA1}. The pull-back of $\chi$ descends to a character of $\pi_1(X_\kappa)$, so our claim follows as in~\cite[prop.~1.3.4(i)]{DeWeilII} by looking at the eigenvalues of the Frobenius operator on the stalks.

\medskip

For $P$ of geometric origin the perverse sheaves $N_\nu\in E(X)$ in part {\em (b)} and hence also all their simple constituents $N$ are of geometric origin. Each such constituent has the form $N \cong L_\varphi \otimes q^*(Q)[\dim(A)]$ by proposition~\ref{prop:chi}(a), so the pullback $i^*(N)$ to $A$ is an isotypic multiple of $i^*(L_\varphi)$ and of geometric origin. Hence $\Pi(i)(\varphi)$ is a torsion character. Writing $\calS(N) = \chi \cdot K(A)$ we can take for $\chi^{-1}$ any torsion character in $\Pi(i)^{-1}(\Pi(i)(\varphi))$. \qed

\medskip

For a homomorphism $f: X\to B$ of abelian varieties, define the relative spectrum $\calS_f(P)$ of a perverse sheaf $P$ on $X$ to be the set of all $\chi \in \Pi(X)$ such that the direct image $Rf_*(K_\chi)$ is not perverse. By abuse of notation, for $\chi \in \Pi(X)$ and $\psi \in \Pi(B)$ we write $\chi \psi  = \chi \cdot (\Pi(f)(\psi)) \in \Pi(X)$. Then the projection formula shows 
\[ Rf_*(P_{\chi \psi}) \;=\; (Rf_*(P_\chi))_\psi, \]
hence $\calS_f(P)$ is invariant under $\Pi(B)$. In particular, if $B=X/A$ for an abelian subvariety $A\subseteq X$, then $\calS_f(P)$ is determined by its image $\overline{\calS}_f(P)$ in $\Pi(A) = \Pi(X) / \Pi(B)$. Furthermore, in theorem~\ref{thm:relative} the assertion {\em for most characters} can be read in $\Pi(A)$, i.e.~theorem~\ref{thm:relative} holds in the stronger sense that $\overline{\calS}_f(P)$ is contained in a finite union of translates of proper algebraic subtori of $\Pi(A)$. Indeed we have 

\medskip \begin{lem} \label{lem:RelativeSpectrum}
$\calS_f(P) \subseteq \calS(P) \cdot \Pi(B)$.
\end{lem} \medskip

{\em Proof.} If $\chi \in \Pi(X)$ does not lie in $\calS(P) \cdot \Pi(B)$, then for any $\psi \in \Pi(B)$ we have $\chi \psi \notin \calS(P)$ and hence $H^\bullet(X, P_{\chi \psi})=H^\bullet(B, (Rf_*(P_\chi))_\psi)$ is not concentrated in degree zero. By theorem~\ref{thm:P} then $Rf_*(P_\chi)$ is not perverse. \qed

\medskip

%
%

\section{Localization at hereditary classes} \label{sec:localization}

In this section we recall certain localization constructions that will be used in what follows to extend our Tannakian results to the case of non-semisimple perverse sheaves. Here $k$ and $\Lambda$ can be arbitrary. For the category $\bfD$ and the functor $rat$ we only require the axiom (D1) from section~\ref{sec:setting}, but we make the following additional assumption:
\begin{itemize}
\item[(T)] {\em Triangulation}. The category $\bfD$ is triangulated and has a $t$-structure with core $\bfP$ which gives rise to the data in (D1).
\end{itemize}
\noindent
We say that a  class $\sfH$ of simple objects in  $\bfP$ is {\it hereditary} if it is stable under the adjoint duality functor $K\mapsto K^\vee$ and if for all~$K\in \sfH$, $L\in \bfD$, $n\in \bbZ$ every simple subquotient of $\pH^n(K*L)$ lies again in $\sfH$. By d\'evissage it suffices of course to check the latter condition only for all simple objects $L\in \bfP$.

\medskip \begin{ex} \label{ex:hereditary}
Suppose that the full subcategory $\bfD^{ss} \subset \bfD$ of all direct sums of degree shifts of simple objects of $\bfP$ satisfies axioms (D1) and (D2) from section~\ref{sec:setting}, and assume for simplicity that $\Lambda$ is algebraically closed. Then the following classes are hereditary:
\begin{enumerate}
 \item[\em (a)] the class~$\sfH_{coh}$ of simple objects $K$ with $H^\bullet(X,K)=0$, \smallskip
 \item[\em (b)] the class $\sfH_{most}$ of simple objects $K$ with $H^\bullet_\chi(X, K)=0$ for most $\chi$, \smallskip \item[\em (c)] the class $\sfH_{Euler}$ of simple objects $K$ with Euler characteristic zero, \smallskip
 \item[\em (d)] the class $\sfH_A$ of simple objects that are invariant under translations by all points in a given abelian subvariety $A\subseteq X$.  \smallskip
\end{enumerate}
\end{ex} \medskip

Indeed, we have $H^\bullet(X,K*L)=H^\bullet(X, K)\otimes H^\bullet(X, L)$ by the K\"unneth formula, and our semisimplicity assumption on the full subcategory $\bfD^{ss}\subseteq \bfD$ ensures that for simple objects $K, L \in \bfP$ the convolution $K*L$ splits into a direct sum of degree shifts of semisimple objects of $\bfP$. So part {\em (a)} is obvious, part {\em (b)} holds because twisting with a character is a tensor functor as we have seen in proposition~\ref{prop:tensorfunctor}, part {\em (c)} follows from corollary~\ref{cor:iso_to_zero} with $\bfD^{ss}$ in place of $\bfD$, and part {\em (d)} is also clear. Note that $\sfH_{most} \subseteq \sfH_{Euler}$ and that for $k=\bbC$ both $\sfH_{coh}$ and $\sfH_{A}$ are contained in $\sfH_{most} = \sfH_{Euler}$ due to corollary \ref{cor:acyclictwist}.

\medskip

For a hereditary class $\sfH$ we denote by $\bfN_\sfH \subseteq \bfD$ the full subcategory of all objects $K\in \bfD$ such that for all $n\in \bbZ$ all simple subquotients of $\pH^n(K)$ are isomorphic to objects in $\sfH$. Then $\bfN_\sfH$ is a thick triangulated tensor ideal of~$(\bfD,*)$, so the localization $\bfD_\sfH = \bfD[\Sigma^{-1}]$ at the class $\Sigma$ of all morphisms with cones in $\bfN_\sfH$ inherits the structure of a rigid symmetric monoidal category such that the localization functor $\bfD \to \bfD_\sfH$ is a tensor functor ACU. 

\medskip

Since $\bfP \cap \bfN_\sfH$  is a Serre subcategory of the abelian category $\bfP$, we can also form the abelian quotient category $\bfP_\sfH = \bfP / \bfP\cap \bfN_\sfH$ as in~\cite[p.~364ff.]{Ga} by inverting all morphisms in $\bfP$ with kernel and cokernel in $\bfP\cap \bfN_\sfH$. 
%
%
The following lemma relates this quotient to the previous localization.

\medskip \begin{lem} \label{rem:gabriel}
The perverse $t$-structure on $\bfD$ induces a $t$-structure on~$\bfD_\sfH$ whose core is the essential image of $\bfP$ under the functor $\bfD\to \bfD_\sfH$, and this essential image is naturally equivalent to the abelian quotient category
$\bfP_\sfH$.
\end{lem} \medskip

For the proof see e.g.~\cite[prop.~3.6.1]{GaL}. In a similar vein we have the following compatibility result for abelian quotient categories.

\medskip \begin{lem} \label{lem:compatible-quotient}
Under the quotient functor $\bfP \to \bfP_{\sfH}$, the image of any Serre subcategory $\bfS \subseteq \bfP$ is naturally equivalent to the abelian quotient category 
\[
 \bfS / \bfS\cap \bfN_\sfH.
\]
\end{lem} \medskip

{\em Proof.} Let us denote by $\bfS_\sfH \subseteq \bfP_\sfH$ the image of $\bfS$. Thus $\bfS_\sfH$ has the same objects as $\bfS$, and by definition of the abelian quotient category $\bfP_\sfH = \bfP/\bfN_\sfH$ the elements of
\[
 \Hom_{\bfS_\sfH}(K, M) \;=\; \Hom_{\bfP_\sfH}(K, M) \quad \textnormal{for} \quad K, M\in \bfS
\]
are given by equivalence classes of diagrams
\[
\xymatrix@R=1.5em@C=3em@M=0.5em{  
 K' \ar[drr] \ar@{^{(}->}[r] & K  \ar@{->>}[r] & N \\
N' \ar@{^{(}->}[r] & M  \ar@{->>}[r] & M'
}
\]
in $\bfP$ with exact rows and with $N, N' \in \bfN_\sfH$. Now by assumption $\bfS \subseteq \bfP$ is a Serre subcategory, so all subquotients of the objects $K, M\in \bfS$ lie again in $\bfS$ and it follows that the above diagram also defines a morphism in the quotient category $\bfS/\bfS\cap \bfN_\sfH$. Furthermore, two diagrams as above are equivalent in~$\bfP$ iff they are equivalent in $\bfS$. Hence the natural functor $\bfS/\bfS\cap \bfN_\sfH \to \bfS_\sfH$ is an equivalence of categories as required. \qed

\section{The Tannaka groups $G(K)$ and $G(X)$} \label{sec:TannakaGroups}

For the definition of the Tannaka groups in theorem~\ref{thm:deligne} we have applied the Andr\'e-Kahn construction to categories of semisimple complexes. It is not clear whether a similar construction works in the non-semisimple case as well, but using theorem~\ref{thm:P} we explain in this section how to define Tannaka groups by another method which also applies to non-semisimple complexes and is compatible with the previous one. Note that non-semisimple perverse sheaves naturally arise as degenerations of semisimple perverse sheaves; in the next section we will provide the appropriate framework to describe such degenerations which in general lead to non-reductive Tannaka groups.


\medskip

Throughout we assume that $\Lambda$ is algebraically closed. Furthermore we require axiom (T) from section~\ref{sec:localization} and the following property:

\begin{itemize}
\item[(S)] {\em Semisimple objects}. The full subcategory $\bfD^{ss} \subseteq \bfD$ of direct sums of degree shifts of simple objects in $\bfP$ satisfies axioms (D1) -- (D3).
\end{itemize}
For instance, by example~\ref{ex:pervgeom} these assumptions are valid for the triangulated category $\bfD = \Dbc(X, \Lambda)$. We want to apply the quotient constructions from the previous section in this axiomatic setting. For the hereditary classes $\sfH = \sfH_\star$ in example \ref{ex:hereditary} with $\star \in \{ \mathit{Euler}, \mathit{coh}\}$ we put
\[
 \bfN_{\star} = \bfN_{\sfH}, \quad \bfD_{\star} = \bfD_\sfH
 \quad \textnormal{and} \quad
 \bfP_{\star} = \bfP_\sfH.
\]
With this notation we have

\medskip \begin{lem}
The category $\bfP_{\mathit{Euler}}$ is a rigid symmetric monoidal abelian subcategory of the rigid symmetric monoidal category $\bfD_{Euler}$.
\end{lem} \medskip

{\em Proof.} The main point is to see that $\bfP_{Euler} * \bfP_{Euler} \subseteq \bfP_{Euler}$. For this it will by d\'evissage suffice to show that the convolution of any two semisimple objects  $K, L\in \bfP$ lies again in $\bfP$ up to direct summands in $\bfN_{Euler}$. But this follows immediately from an application of theorem~\ref{thm:mult} and corollary~\ref{cor:iso_to_zero} to the category $\bfD^{ss}$ in place of $\bfD$, using our assumption (S).
\qed

\medskip

We can now generalize our earlier construction of Tannakian categories in the following way. Let $\bfP^\chi \subset \bfP$ be the full abelian subcategory consisting of all objects $P$ such that all simple subquotients $Q$ of $P$ satisfy $H^i_\chi(X, Q)=0$ for all cohomology degrees $i\neq 0$. We know from lemma~\ref{lem:compatible-quotient} that the abelian quotient category
\[
 \bfP^\chi_{coh} \;:=\; \bfP^\chi / \bfP^\chi \cap \bfN_{coh}\;=\; \bfP^\chi / \bfP^\chi \cap \bfN_{Euler} 
\]
is naturally equivalent to the image of $\bfP^\chi$ inside $\bfP_{Euler}$.

\medskip \begin{thm}
The category $\bfP_{coh}^\chi$ is a rigid symmetric monoidal abelian subcategory of the rigid symmetric monoidal category $\bfP_{coh}$. Furthermore, we have an equivalence
\[
 \bfP_{coh}^\chi \; \cong \; \Rep_\Lambda(G)
\]
with the rigid symmetric monoidal abelian category of finite-dimensional linear representations of some affine group scheme $G=G(X, \chi)$ over $\Lambda$.
\end{thm} \medskip

{\em Proof.} To see that $\bfP_{coh}^\chi * \bfP_{coh}^\chi \subseteq \bfP_{coh}^\chi$, it will by d\'evissage suffice to see that the convolution of any two semisimple objects of $\bfP^\chi$ lies again in $\bfP^\chi$ up to a direct summand in $\bfN_{coh}$. But this follows from proposition~\ref{prop:tensorfunctor} and from the K\"unneth formula, using the semisimplicity axiom (S). So $\bfP_{coh}^\chi$ is a rigid symmetric monoidal abelian subcategory of $\bfP_{Euler}$. For the remaining statement it suffices by~\cite[th.~2.11]{DM} to find a fibre functor, i.e.~an exact, faithful, $\Lambda$-linear tensor functor from $\bfP_{coh}^\chi$ to the category of finite-dimensional vector spaces over $\Lambda$. But for this we can take the functor $K\mapsto H^0_\chi(X, K)$ which is well-defined on $\bfP^\chi_{coh}$ because it vanishes on $\bfP^\chi \cap \bfN_{coh}$.
\qed

%
%




\medskip

In the analytic case where $k=\Lambda = \bbC$, the vanishing theorem~\ref{thm:P} shows that every object $K\in \bfP$ is contained in a subcategory of the form $\bfP^\chi$ for some character $\chi$ of the fundamental group. In any case, we have

\medskip \begin{cor}
For $K\in \bfP^\chi$ let $\langle K \rangle \subseteq \bfP^\chi$ denote the rigid symmetric monoidal abelian subcategory generated in the quotient category $\bfP_{coh}^\chi$ by the subquotients of $(K\oplus K^\vee)^{*r}$ with $r\in \bbN_0$. Then we have an equivalence
\[
 \langle K \rangle \; \cong \; \Rep_\Lambda(G)
\]
with the rigid symmetric monoidal abelian category of finite-dimensional linear representations of some affine algebraic group $G=G(K)$ over $\Lambda$.
\end{cor} \medskip

{\em Proof.} This follows directly from the previous theorem via the Tannakian formalism~\cite[prop.~2.20(b)]{DM}. In fact the affine algebraic group $G(K)$ is a quotient of the affine group scheme $G(X, \chi)$.
\qed

\medskip

The notation in the above corollary is slightly ambiguous: At least in the case where~$k=\Lambda=\bbC$, theorem~\ref{thm:P} implies that each object $K\in \bfP$ lies in $\bfP^\chi$ for many different characters $\chi$. However, since we assumed the coefficient field~$\Lambda$ to be algebraically closed, we have the following

\medskip \begin{rem} \label{lem:change_of_fiber_functor}
Up to isomorphism, the algebraic group $G(K)$ depends only on the object $K$ but not on the character $\chi$.
\end{rem} \medskip

{\em Proof.} The only reason why we have defined $\langle K \rangle$ as a subcategory of~$\bfP^\chi_{coh}$ is that we wanted to have a fibre functor. Indeed we have already remarked above that $\bfP^\chi_{coh} \subseteq \bfP_{Euler}$, so as an abstract $\Lambda$-linear rigid symmetric monoidal abelian category the category $\langle K \rangle$ does not depend on the character $\chi$ that we have chosen. In other words, the choice of the character only affects the fibre functor on our category. But since $G=G(K)$ is an affine algebraic group over the algebraically closed field $\Lambda$, any two fibre functors for the Tannakian category $\Rep_\Lambda(G)$ are isomorphic by~\cite[th.~3.2]{DM}. \qed

%
%

\medskip 

%
%

We also remark that if the perverse sheaf $rat(K)_\chi \in \Perv(X, \Lambda)$ is defined over a subfield $\Lambda_0 \subseteq \Lambda$, then our fibre functor descends to $\Lambda_0$. This allows to define the Tannaka group of $\langle K \rangle$ as an affine algebraic group over~$\Lambda_0$, but as such its isomorphism class may depend on~$\chi$.

\medskip \begin{lem}
For $K \in \bfD^{ss}\cap \bfP$  the Tannaka group $G(K)$ is isomorphic over $\Lambda$ to the Tannaka group $G(\bfT)$ of the
symmetric monoidal subcategory $\bfT$ generated by $K$ in the André-Kahn quotient
of $\bfD^{ss}$ as in section \ref{sec:tannakian}.
\end{lem} \medskip

{\it Proof}. Let $\bfP^{ss}_{coh}$ be the essential image of $\bfP^{ss}=\bfP\cap \bfD^{ss}$ in $\bfP_{coh}$, and denote by $\bfPbar^{ss}$ the image of $ \bfP^{ss}$ in the Andr\'e-Kahn quotient of $\bfD^{ss}$ as defined in section~\ref{sec:andre-kahn}. Then $\bfPbar^{ss}$ is a semisimple abelian category, and also a symmetric monoidal category by theorem~\ref{thm:mult}. By corollary~\ref{cor:exact} the functor $\bfP^{ss} \to \bfPbar^{ss}$ is exact. So if $s$ is a morphism in $\bfP^{ss}$ whose kernel and cokernel lie in~$\bfN_{coh}$, then the kernel and cokernel of the corresponding morphism in $\bfPbar^{ss}$ are zero, i.e.~the morphism $s$ becomes invertible in $\bfPbar^{ss}$. From the universal property of the localization $\bfP^{ss}_{coh}$ we thus obtain a unique functor $p$ such that the following diagram commutes, where ${q}$ and $\overline{q}$ denote the natural quotient functors.
\[
\xymatrix@C=1.5em@R=1.5em{ &   \bfP^{ss} \ar[dr]^{\overline q} \ar[dl]_{q} & \\  \bfP^{ss}_{coh} \ar@{..>}[rr]_{\exists ! \, p}  & 
 & \bfPbar^{ss} }
\]
Now $p$ is a tensor functor ACU since ${q}$ and $\overline{q}$ are,  and it induces a functor $p_K$ from the semisimple abelian subcategory~$\langle K \rangle \subseteq  \bfP^{ss}_{coh}$ to $\bfT \subseteq \bfPbar^{ss}$ which is essentially surjective since it is the identity on objects. One easily checks that this functor is fully faithful, hence an equivalence of categories.
\qed

\medskip

For the rest of this section we will be concerned only with the special case where $\bfD = \Dbc(X, \Lambda)$. In this case we consider the pro-algebraic group
\[ G(X) \;=\; G(X, 1) \]
attached to the trivial character $\chi = 1$. The reason why we only consider the trivial character is of course that the groups $G(X, \chi)$ are isomorphic for all choices of $\chi$, indeed for $\bfP = \Perv(X, \Lambda)$ the categories $\bfP_{coh}^1$ and $\bfP_{coh}^\chi$ are isomorphic to each other via the twisting functor $K\mapsto K_{\chi^{-1}}$.

%

\medskip \begin{lem} Suppose $k=\Lambda =\bbC$. Then the maximal abelian quotient group 
of the group of connected components of $G(X)$ is 
$$ \pi_0(G(X))^{ab} \ =\ \pi_1^{et}(\hat X,0)(-1)$$
where $\hat X$ is the dual abelian variety of $X$. Here the Tate twist~$(-1)$ refers to the action of $\mathrm{Gal}(k/k_0)$ for any subfield $k_0 \subseteq k$ over which $X$ is defined.
\end{lem} \medskip

{\em Proof.}
For $K\in \bfP^1$ the epimorphism $G(K) \twoheadrightarrow \pi_0(G(K))^{ab}$ defines the full subcategory 
$\Rep_\Lambda(\pi_0(G(K))^{ab}) \subseteq \Rep_\Lambda(G(K)) = \langle K \rangle$
generated by the characters of $G(K)$ of finite order. For $\bfD = \Dbc(X, \Lambda)$ and $k=\Lambda=\bbC$, any character of $G(K)$ is by part~{\em (b)} of proposition \ref{prop:chi} represented by a skyscraper sheaf $\delta_x$ supported in a point $x\in X(\bbC)$. Since $\delta_x*\delta_y =\delta_{x+y}$, such a character has finite order iff $x$ is a torsion point in $X(\bbC)$. Hence 
\[ \pi_0(G(X))^{ab}(1) \;=\;
\varprojlim_{n} \bigl( \Hom(X(\bbC)[n],\bbG_m) \bigr) \;=\; \pi_1^{ et}(\hat X,0) \]
by Pontryagin duality, and we are done.
\qed

\medskip \begin{lem} \label{lem:homomorphisms}
Every homomorphism $f:X\longrightarrow Y$ of abelian varieties over~$k$ induces a homomorphism of pro-algebraic groups $$G(f): \quad G(Y) \longrightarrow G(X).$$
If $f$ is surjective, then this homomorphism $G(f)$ is a closed embedding. 
\end{lem} \medskip

{\em Proof.}
One easily checks $Rf_*(K*L)=Rf_*(K)*Rf_*(L)$. Furthermore we have $Rf_*(\bfN_{coh}) \subseteq \bfN_{coh}$. Indeed, by d\'evissage it suffices to check this property for simple perverse sheaves, where it follows from the decomposition theorem. One then deduces that $Rf_*$ induces a tensor functor~ACU from $\Perv(X, \Lambda)^1_{coh}$ to $\Perv(Y, \Lambda)^1_{coh}$ and hence a homomorphism $G(f)$.

\medskip

If $f$ is surjective, $G(f)$ is a closed embedding. Indeed, by \cite[prop.~2.21b)]{DM}  it is enough to show that for perverse $K$ on $Y$ there exists a perverse sheaf~$L$ on $X$ such that $K$ is a retract of $Rf_*(L)$. This assertion can be reduced  to the cases where $f$ is either an isogeny or a projection $X=Z\times Y \to Y$ onto a factor. In these two cases one can take $L=f^*(K)$ resp.~$L={\bf1} \boxtimes K$. \qed

\section{Nearby Cycles} \label{sec:nearby}

In this section we describe how the Tannaka groups $G(K)$ vary in algebraic families. This paves the way for degeneration arguments and is so far the most efficient tool to obtain information about the arising Tannaka groups in concrete geometric situations, see for example~\cite{KrW}. For $i\in \{0, 1\}$, let $X_i$ be an abelian variety over an algebraically closed field $k_i$ which has characteristic zero or is the algebraic closure of a finite field. Let $\Lambda = \bbC$ or $\Lambda = \Qbar_l$ for some prime number $l\neq \mathit{char}(k_0), \mathit{char}(k_1)$. Suppose we have $\Lambda$-linear rigid symmetric monoidal categories $\bfD_i$ 
and faithful $\Lambda$-linear tensor functors ACU $rat_i: \bfD_i \to \Dbc(X_i, \Lambda)$, and assume that both categories $\bfD_0$ and $\bfD_1$ satisfy the axioms (T) and (S) from the previous sections.

\medskip

Let $\psi_\bfD: \bfD_0 \longrightarrow \bfD_1$ and $\psi_D: \Dbc(X_0, \Lambda) \longrightarrow \Dbc(X_1, \Lambda)$ be triangulated and $t$-exact tensor functors ACU such that the diagram
\[
\xymatrix@C=4em@M=0.5em{
 \bfD_0 \ar[d]_{\psi_\bfD} \ar[r]^-{rat} & \Dbc(X_0, \Lambda) \ar[d]^{\psi_D} \\
 \bfD_1 \ar[r]^-{rat} & \Dbc(X_1, \Lambda)
}
\]
commutes, and suppose we have functorial isomorphisms 
\[ H^\bullet(X_1, \psi_D(K)) \; \cong \; H^\bullet(X_0, K) \quad \textnormal{for all} \quad K\in \Dbc(X_0, \Lambda). \]
Assume furthermore that we have an identification $\pi_1(X_0, 0) = \pi_1(X_1, 0)$ under which $\psi_D((-)_\chi) = (\psi_D(-))_\chi$ for all characters $\chi$ of this group.   In what follows we simply write $\psi$ for both $\psi_\bfD$ and~$\psi_D$.

\medskip \begin{lem} \label{lem:psi}
Let $K \in \bfP_0$. If $\psi(K) \in \bfP_1^\chi$ for some character $\chi$, then we have a closed embedding (depending on the character)
\[
 G(\psi(K)) \; \hookrightarrow \; G(K).
\] 
\end{lem} \medskip

{\em Proof.} Since the functor $\psi$ is exact and compatible with hypercohomology and character twists in the sense explained above, it restricts to a functor from~$\bfP_0^\chi$ to $\bfP_1^\chi$ which sends $\bfP_0^\chi \cap \bfN_{coh}$ into $\bfP_1^\chi\cap \bfN_{coh}$. Hence $\psi$ induces a tensor functor ACU between the corresponding quotient categories and in particular between their subcategories $\langle K \rangle$ and $\langle \psi(K) \rangle$. So we can proceed as in the proof of lemma~\ref{lem:homomorphisms}.
\qed

\medskip

The above result applies in particular in the following situation. Let $S$ be the spectrum of a Henselian discrete valuation ring with closed point~$s$ and generic point~$\eta$. Let $\etabar$ be a geometric point over $\eta$. Let $\Sbar$ be the normalization of $S$ in the residue field $\kappa(\etabar)$, and let $\sbar$ be a geometric point of $\Sbar$ over $s$. For an abelian scheme $X$ over $S$, put $\Xbar = X\times_S \Sbar$. Consider the commutative diagramm%
\[
 \xymatrix@C=2.5em@M=0.5em{
 X_\sbar \ar@{^{(}->}[r]^\ibar \ar[d] & \Xbar \ar[d] & X_\etabar \ar@{_{(}->}[l]_\jbar \ar[d] \\
 \sbar \ar@{^{(}->}[r] & \Sbar & \etabar \ar@{_{(}->}[l] 
}
\]
where $\ibar$ and $\jbar$ are the natural morphisms. We then have the functor of nearby cycles \cite[exp. XIII-XIV]{SGA7}
\[
 \psi \;=\; \ibar^* R\jbar_*: \quad \bfD_0 \;=\; \Dbc(X_\etabar, \Lambda) \;\longrightarrow \; \bfD_1 \;=\; \Dbc(X_\sbar, \Lambda).
\]
This functor is $t$-exact for the perverse $t$-structures by~\cite[cor.~4.5]{Il}, and theorem~4.7 in loc.~cit.~implies that it is a tensor functor for convolution. Note that $H^\bullet(X_\sbar, \psi (K)) = H^\bullet(X_\etabar, K)$, for all $K\in \Dbc(X_\etabar, \Lambda)$. Furthermore, since $\Xbar$ is proper and smooth over $\Sbar$, by~\cite[exp.~X, cor.~3.9]{SGA1} we have a specialization epimorphism 
$
 \mathit{sp}: \pi_1(X_\etabar, 0) \twoheadrightarrow \pi_1(X_\sbar, 0)
$
whose kernel is  a pro-$p$-group for $p=\mathit{char}(\kappa(\sbar))$. If $\kappa(\sbar)$ has characteristic zero, then $\mathit{sp}$ is an isomorphism. Extending local systems on $X_\etabar$ to local systems on $\Xbar$, one then sees that $\psi(L_\chi)=L_{\chi \circ \mathit{sp}^{-1}}$ for any character $\chi$ of $\pi_1(X_\etabar, 0)$. In this case we also write~$\chi$ for the character $\chi \circ \mathit{sp}^{-1}$ of $\pi_1(X_\sbar, 0)$ by abuse of notation, so that $\psi(K_\chi)=(\psi(K))_\chi$ for all characters $\chi$ and all $K\in \Dbc(X_\etabar, \Lambda)$.

\section{Appendix: Reductive supergroups} \label{sec:supergroups}

In this appendix we recall the definition of an algebraic super group and collect some basic facts about these in the reductive case. Throughout let $\Lambda$ be an algebraically closed field of characteristic zero. As in~\cite[p.~16]{WeSS} we consider triples $\bfG = (G, \frakg_-, Q)$ consisting of 
\begin{itemize}
 \item a classical algebraic group $G$ over $\Lambda$, whose Lie algebra equipped with the adjoint action of $G$ we denote by $\frakg_+=\Lie(G)$,
 \item a finite-dimensional algebraic representation $\frakg_-$ of $G$ over $\Lambda$, given by a homomorphism $Ad_-: G\to \Gl(\frakg_-)$,
 \item a $G$-equivariant quadratic form $Q: \frakg_- \to \frakg_+, \, Q(v)= [v,v]$ defined by a symmetric $\Lambda$-bilinear form $[\cdot, \cdot ]: \frakg_-\times \frakg_- \to \frakg_+$.
\end{itemize}
Such a triple $\bfG$ is called an algebraic super group over $\Lambda$, if the differential $ad_-=\Lie(Ad_-)$ of $Ad_-$ satisfies $ad_-(Q(v)) (v) = 0$ for all~$v\in \frakg_-$. We define a homomorphism 
\[ h: \; (G_1, \frakg_{1-}, Q_1) \longrightarrow (G_2, \frakg_{2-}, Q_2) \] 
of algebraic super groups over $\Lambda$ to be a pair $h=(f, g)$, where $f: G_1\to G_2$ is a homomorphism of algebraic groups and $g: \frakg_{1-}\to \frakg_{2-}$ a $\Lambda$-linear and \mbox{$f$-equivariant} map such that $Q_2 \circ g = \Lie(f) \circ Q_1$. Such a homomorphism is a mono- resp.~an epimorphism of algebraic super groups iff both $f$ and $g$ are mono- resp.~epimorphisms. We define the parity automorphism of an algebraic super group $\bfG = (G, \frakg_-, Q)$ to be $h=(\id_{\, G}, -\id_{\, \frakg_-}): \bfG\to \bfG$. These constructions are motivated by the following example.

\medskip

Let $A=A_+\oplus A_-$ be an affine super Hopf algebra over~$\Lambda$, i.e.~a graded commutative $\bbZ/2\bbZ$-graded Hopf algebra of finite type over~$\Lambda$. Let $J\trianglelefteq A$ be the ideal generated by $A_-$. Then $G=Spec \, A/J$ is a classical algebraic group, and the left invariant super derivations of $A$ form a super Lie algebra $\frakg = \frakg_+ \oplus \frakg_-$ with a natural action of $G$ extending the adjoint action on $\frakg_+=\Lie(G)$.
If we take $Q(v)=[v,v]$ for the super bracket $[\cdot, \cdot]: \frakg \times \frakg \to \frakg$, then by loc.~cit.~$\bfG = (G, \frakg_-, Q)$ is an algebraic super group over $\Lambda$. By 
loc. cit. this realizes the opposite of the category of affine super Hopf algebras as a full subcategory of the category of algebraic super groups. Hence for algebraic super groups associated to affine super Hopf algebras, the notions introduced here are compatible with those in~\cite{DelCT}.

\medskip

A particular instance is the general linear super group $\bfG = \bfGl(V)$ attached to a super vector space $V=V_+\oplus V_-$ of finite dimension over~$\Lambda$. In this case $G=\Gl(V_+) \times \Gl(V_-)$, $\frakg_- = \Hom_\Lambda(V_+, V_-) \oplus \Hom_\Lambda(V_-, V_+) \subset \End_\Lambda(V)$ with the adjoint action of $G$, and one takes $Q(A\oplus B) = AB + BA$. 

\medskip

For any algebraic super group $\bfG = (G, \frakg_-, Q)$ over $\Lambda$, let~$\bfG^0 = (G^0, \frakg_-, Q)$ denote its Zariski connected component, and define its super center to be $Z(\bfG) = (Z, 0, 0)$ where $Z\subseteq Z(G)$ is the largest central subgroup of $G$ acting trivially on $\frakg_-$. For $g\in G$ we put $\mathrm{int}(g) = (g^{-1}\, (-) \, g, Ad_-(g)): \bfG \to \bfG$. Then $Z\subseteq G$ is the subgroup of all $g\in G$ such that $\mathrm{int}(g)=\id_\bfG$.

\medskip

A super representation of $\bfG$ is a homomorphism $\rho_V: \bfG \to \bfGl(V)$ for some super vector space $V$. By definition, a homomorphism between two super representations $\rho_V$ and $\rho_W$ is a homomorphism $V\to W$ of super vector spaces such that the induced homomorphism $h: \bfGl(V)\to \bfGl(W)$ satisfies $h\circ \rho_V = \rho_W$. The category $\Rep_\Lambda(\bfG)$ of super representations of $\bfG$ over $\Lambda$ is an abelian $\Lambda$-linear rigid symmetric monoidal category with respect to the super tensor product. We also have Schur's lemma:

\medskip \begin{lem}
Let $\rho_V: \bfG \to \bfGl(V)$ be an irreducible super representation. Then every endomorphism $\varphi$ of $\rho_V$ has the form $\varphi = \lambda\cdot \id_V$ for some $\lambda \in \Lambda$. 
\end{lem} \medskip

{\em Proof.} The proof works as in the classical case. Notice that by definition we only consider endomorphisms preserving the~$\bbZ/2\bbZ$-grading. Otherwise Schur's lemma would have to be modified, see~\cite[prop.~2, p.~46]{Sc}. 
\qed

\medskip

In particular, it follows that the super center $Z(\bfG)$ acts on any irreducible super representation of $\bfG$ by a character $\chi: Z(\bfG) \to \Lambda^*$ (recall that the super center is a classical commutative algebraic group and that each of its elements defines an endomorphism of any given super representation of~$\bfG$).

\medskip

An algebraic super group $\bfG = (G, \frakg_-, Q)$ is called {\it reductive}, if the abelian category $\Rep_\Lambda(\bfG)$ is semisimple. Let us briefly recall the classification of reductive super groups from~\cite{WeSS}. Every classical reductive group $G$, viewed as a super group with $\frakg_-=0$, is a reductive super group. Other examples include the orthosymplectic super groups ${\bf Spo}_\Lambda(2r, 1)$ with $r\in \bbN$, defined as follows: Fix a non-degenerate antisymmetric $2r\times 2r$ matrix $J$ over $\Lambda$, and consider $\Sp_\Lambda(2r, J) = \{g\in \Gl_\Lambda(2r) \mid g^t J g = J\}$. Then
\[
 {\bf Spo}_\Lambda(2r, 1) \;=\; (\Sp_\Lambda(2r, J), \Lambda^{2r}, Q) \ ,
\]
where $\Lambda^{2r}$ is equipped with the standard action of $\Sp_\Lambda(2r, J)$ and where the map $Q: \Lambda^{2r} \to \Sp_\Lambda(2r, J)$ is defined by $Q(v)_{ik} = \sum_{j=1}^{2r} v_i v_j J_{jk}$. A different choice of the matrix $J$ gives an isomorphic super group.

\medskip

In general, by theorem~6 of loc.~cit., an algebraic super group $\bfG$ over~$\Lambda$ is reductive iff there exists a classical reductive group $H$ and $N\in \bbN_0$, $n_i, r_i\in \bbN$ such that $\bfG$ is isomorphic to a semidirect product
\[
 \bfG \;=\; \Bigl( \prod_{i=1}^N \bigl( {\bf Spo}_\Lambda(2r_i, 1) \bigr)^{n_i} \Bigr) \rtimes H
\]
defined by a homomorphism $\pi_0(H) \longrightarrow \prod_{i=1}^N \mathfrak{S}_{n_i}$ where each symmetric group~$\mathfrak{S}_{n_i}$ acts on $\bigl({\bf Spo}_\Lambda(2r_i, 1) \bigr)^{n_i}$ by permutation of the factors. 

\medskip \begin{cor} \label{cor:reductive}
For any reductive super group $\bfG$ over $\Lambda$, the underlying classical group $G$ is reductive, and the super center $Z(\bfG)$ is a subgroup of finite index in the center $Z(G)$.
\end{cor} \medskip

{\em Proof.} By the above it suffices to show this in the case $\bfG = {\bf Spo}_\Lambda(2r, 1)$ for some $r\in \bbN$. But then $G=\Sp_\Lambda(2r)$, and $Z(G)=\mu_2$ is finite. \qed

\medskip

\medskip \begin{prop} \label{prop:toruslift}
Let $h: \bfG_1 \to \bfG_0$ be a homomorphism of reductive super groups over $\Lambda$ which induces an epimorphism $f: G_1\twoheadrightarrow G_0$ on the underlying classical groups. If the super center $Z(\bfG_0)$ contains a classical torus $T_0$, then $Z(\bfG_1)$ contains a classical torus $T_1$ such that $p$ induces an isogeny~$p: T_1 \to T_0$.
\end{prop} \medskip

{\em Proof.} The category of tori (or diagonalizable groups) over $\Lambda$, up to isogeny, is
equivalent to the category of finite-dimensional vector spaces over $\bbQ$ via the  cocharacter functor 
$T \mapsto X(T)=Hom(\bbG_m,T)\otimes_\bbZ \bbQ\ .$ 
If $\pi$ is a finite group acting on $T$, then $X((T^\pi)^0) = X(T)^\pi$  for the fixed group $T^\pi$.

\medskip
 
For reductive super groups $\bfG$ note 
$Z(\bfG)^0 = Z(G)^0$ by corollary~\ref{cor:reductive}. On $Z(G^0)$ the group $G$ acts by conjugation, thus defines an
action of the finite group $\pi=\pi_0(G)$. By definition $Z(G)^0 \subset Z(G^0)^\pi$
and this is a subgroup of finite index: $Z(G)^0 = (Z(G^0)^\pi)^0$. This follows by an application of the cocharacter functor since $X(Z(G)^0)
= X(Z(G^0))^\pi$.

\medskip
For the proof of the proposition it suffices to show that $T_0$ is contained in the image of $Z(G_1)^0$.
By assumption $h$ induces an epimorphism $f:G_1\twoheadrightarrow G_0$ of classical reductive groups, hence an epimorphism $(G_1)^0 \twoheadrightarrow (G_0)^0$ of their connected components.  
By the theory of classical connected reductive groups the torus $T=Z((G_0)^0)^0$ is the image of the torus $S=Z((G_1)^0)^0$. The epimorphism $f:S \twoheadrightarrow T$ is equivariant for the action of $\pi=\pi_0(G_1)$ on $S$ and~$T$, where the latter is induced by the homomorphism $\pi_0(G_1)\to \pi_0(G_0)$. We claim that 
$  h: (S^{\pi})^0 \to (T^\pi)^0 $
is surjective. Indeed, the functor
of invariants under a finite group $\pi$ is right exact on the category of finite-dimensional vector spaces over $\bbQ$. Since $(S^{\pi})^0 \subset Z(G_1)^0$ and $T_0 \subset (T^\pi)^0$
this completes the proof. \qed

\end{document}